\documentclass{article}

\usepackage{epsfig}
\usepackage{amsfonts}
\usepackage{cite}
\usepackage{epstopdf}
\usepackage{mathptmx}

\newtheorem{theorem}{Theorem}[section]
\newtheorem{example}[theorem]{Example}
\newtheorem{remark}[theorem]{Remark}

\begin{document}

\title{Traveling wave solutions of some important Wick-type 
fractional stochastic nonlinear partial differential equations\thanks{This is a preprint 
of a paper whose final and definite form is with 'Chaos, Solitons \& Fractals', 
ISSN 0960-0779 [\texttt{https://doi.org/10.1016/j.chaos.2019.109542}]. 
Submitted 19-Sept-2019; Revised 14-Nov-2019; Accepted 18-Nov-2019.}}

\author{Hyunsoo Kim$^{1,}$\thanks{E-mail Address: hskiminkorea@gmail.com; ORCID: 0000-0003-0819-7496}\ ,
Rathinasamy Sakthivel$^{2,}$\thanks{E-mail Address: krsakthivel@buc.edu.in; ORCID: 0000-0002-5528-2709}\ ,
Amar Debbouche$^{3,4,}$\thanks{Corresponding author; E-mail Address: amar\_debbouche@yahoo.fr; ORCID: 0000-0003-4321-9515}\ ,
Delfim F. M. Torres$^{4,}$\thanks{E-mail Address: delfim@ua.pt; ORCID: 0000-0001-8641-2505}}

\date{{\small $^{1}$Department of Mathematics,
Sungkyunkwan University, Suwon 16419, Republic of Korea}\\
{\small $^{2}$Department of Applied Mathematics,
Bharathiar University, Coimbatore 641 046, India}\\
{\small $^{3}$Department of Mathematics, Guelma University, Guelma 24000, Algeria}\\ 
{\small $^{4}$CIDMA -- Center for Research and Development in Mathematics and Applications,}\\
{\small Department of Mathematics, University of Aveiro, Aveiro 3810-193, Portugal}}

\maketitle


\begin{abstract}
In this article, exact traveling wave solutions of a Wick-type stochastic nonlinear
Schr\"{o}dinger equation and of a Wick-type stochastic fractional Regularized 
Long Wave-Burgers (RLW-Burgers) equation have been obtained 
by using an improved computational method. 
Specifically, the Hermite transform is employed for transforming Wick-type 
stochastic nonlinear partial differential equations into deterministic nonlinear
partial differential equations with integral and fraction order. Furthermore, 
the required set of stochastic solutions in the white noise space is obtained  
by using the inverse Hermite transform. Based on the derived solutions,
the dynamics of the considered equations are performed with some
particular values of the physical parameters. The results reveal that the  
proposed improved computational technique can be applied to solve various 
kinds of Wick-type stochastic fractional partial differential equations.
\end{abstract}



\section{Introduction}

Stochastic nonlinear partial differential equations
play an important role in various fields of science and engineering
\cite{pang,kafash,lee}. In general, finding solutions 
of stochastic nonlinear partial differential equations is
more complex due to its additional random terms 
when compared to deterministic partial
differential equations \cite{book1,book,300,ref1}. 
In particular, the exact traveling wave solution 
receives a great deal of attention since it reveals the internal
mechanism of physical phenomena and also it can be used to assist
numerical solvers to compare the correctness of their results.
In particular, the study on obtaining exact traveling wave solutions 
of stochastic nonlinear partial differential equations via computational 
techniques is still in the initial stage. At present, only few works 
have been reported regarding the exact traveling wave solutions 
of stochastic nonlinear partial differential equations 
\cite{kim2,kim3,w3,w4,w5,w6}. Saha Ray and Singh
obtained new exact solutions of the Wick-type stochastic
Kudryashov--Sinelshchikov equation with the aid of improved
Sub-equation method \cite{w2}. In addition, exact traveling wave 
solutions for the Wick-type stochastic coupled KdV-mKdV equations 
have been obtained by employing the Jacobian elliptic 
function expansion method in \cite{w1}.

Consider the variable coefficient nonlinear 
Schr\"{o}dinger equation in the following form:
\begin{equation}
\label{pertur-schrodinger}
i \psi_t +\alpha(t) \psi_{xx} + \beta(t)\psi |\psi|^{2} + \lambda(t)
\psi =0, \quad i^2 = -1,
\end{equation}
which is used to describe a variety of significant phenomena in ocean dynamics, 
fluid mechanics and plasma physics. Here,
$(x,t) \in \mathbb{R} \times \mathbb{R}^{+}$; $\alpha(t),
\beta(t)$ and $\lambda(t)$ are bounded measurable or integrable
functions on $\mathbb{R}^{+}$ \cite{yazid,gung}.
Based on Eq. (\ref{pertur-schrodinger}), we consider the
Wick-type stochastic nonlinear Shr\"{o}dinger equation in form
\begin{equation}
\label{stochatic-schrodinger}
i \Psi_t +\alpha^*(t) \diamondsuit \Psi_{xx} + \beta^*(t)
\diamondsuit \Psi \diamondsuit |\Psi|^{\diamondsuit 2} + \lambda^*(t)
\diamondsuit \Psi =0,
\end{equation}
where $\diamondsuit$ is the Wick product on the Kondratiev
distribution space $(S(\mathbb{R}^d))^*$ and $\alpha^*(t)$, $\beta^*(t)$ and
$\lambda^*(t)$ are white noise functionals defined in $(S(\mathbb{R}^d))^*$, 
where $(S(\mathbb{R}^d))^*$ is the white noise functional space \cite{tawil,dai2,pan,yin}.

On the other hand, fractional calculus is devoted to the study of integrals 
and derivatives of fractional order, which has been used to describe nonlinear 
physical models arising in real world problems such as in viscoelastic systems, 
complex biological systems, signal processing control and image processing 
of computer vision \cite{moha,tate,garra,capo,ostal,povs,mitter}. Fractional 
nonlinear partial differential equations (Fractional NPDEs) can be written as 
generalizations of classical NPDEs through several senses of fractional derivatives 
\cite{ref2,ref3,ref4,ref5,ref6}, including the recent popular one called the 
Atangana--Baleanu fractional derivative in Caputo sense, which was introduced in \cite{ABC}. 
It should be noted that fractional NPDEs can characterize the physical phenomena 
better than deterministic NPDEs \cite{fra2,zhang1}. In particular, it should be noted 
that it is more difficult to solve stochastic NPDEs in comparison 
with deterministic NPDEs. Though exact traveling wave solutions of
Wick-type stochastic NPDEs have been extensively studied by many researchers, 
finding exact traveling wave solutions of fractional versions of them is still in
the initial stage. Motivated by this consideration, in this paper we consider the fractional 
Regularized Long Wave-Burgers equation (RLW-Burgers) with Burgers-type dissipative 
and dispersive terms in the form
\begin{equation}
\label{rlw-frac}
\frac{\partial^{\alpha} V}{\partial t^{\alpha}} +
p(t)\frac{\partial^{\alpha} V}{\partial x^{\alpha}}
+q(t)V\frac{\partial^{\alpha} V}{\partial x^{\alpha}} +r(t)
\frac{\partial^{2\alpha} V}{\partial
x^{2\alpha}}+s(t)\frac{\partial^{3\alpha} V}{\partial t^{\alpha}
\partial x^{2\alpha}} =0, 
\end{equation}
where $p(t)$, $q(t)$, $r(t)$ and $s(t)$ are bounded measurable functions 
on $\mathbb{R}^+$ and the fractional order $\alpha$ 
belongs to $(0,1)$ \cite{zhang1,bona}.

Motivated by Eq. (\ref{rlw-frac}), in this paper we consider 
the Wick-type fractional stochastic version of Eq.
(\ref{rlw-frac}) in the following form:
\begin{eqnarray}
\frac{\partial^{\diamondsuit\alpha} V}{\partial t^{\diamondsuit\alpha}} 
+ P(t) \diamondsuit \frac{\partial^{\diamondsuit\alpha} V}{\partial
x^{\diamondsuit\alpha}} +Q(t) \diamondsuit V \diamondsuit
\frac{\partial^{\diamondsuit\alpha} V}{\partial x^{\diamondsuit\alpha}}
+R(t) \diamondsuit \frac{\partial^{\diamondsuit 2\alpha} V}{\partial
x^{\diamondsuit 2\alpha}}+S(t) \diamondsuit
\frac{\partial^{\diamondsuit 3\alpha} V}{\partial t^{\diamondsuit\alpha}
\partial x^{\diamondsuit 2\alpha}} =0, 
\label{rlw-wick}
\end{eqnarray}
where ``$\diamondsuit$'' is the Wick product on the Kondratiev
distribution space $(S(\mathbb{R}^d))^*$, and $P(t), Q(t), R(t)$ and
$S(t)$ are white noise functionals defined in $(S(\mathbb{R}^d))^*$ \cite{holden}. 
Also, $P(t)=f_1(t)+c_1 W(t)$, $Q(t)=f_2(t)+c_2 W(t)$, $R(t)=f_3(t)+c_3 W(t)$,
$S(t)=f_4(t)+c_4 W(t)$ and $\left\{f_i(t) \right\}_{i=1}^4$ are bounded
and integrable functions on $\mathbb{R}^{+}$,  where $\left\{c_i \right\}_{i=1}^4$ denote
constants and $W(t)$ is the Gaussian white noise that satisfies $W(t) =
\dot{B}_t = d B(t)/dt$ where $B(t)$ is a Brownian motion.

The main aim of this work is to obtain new Wick-type exact traveling wave solutions
for Eq. (\ref{stochatic-schrodinger}) and Eq. (\ref{rlw-wick}) by using the 
improved computational algorithm and the modified fractional
sub-equation method. Further, we obtain the white noise functional
solutions of Eq. (\ref{rlw-wick}) based on the Hermite transform and
the Caputo fractional derivative. 
This paper is organized as follows. In Section~\ref{sec:2}, we
describe the Wick-type product and the improved computational
algorithm and we obtain new Wick-type exact traveling wave solutions 
of the Wick-type stochastic nonlinear Schr\"{o}dinger equations, 
and we perform the dynamics of the given solutions of the considered equation. 
In Section~\ref{sec:3}, we describe the modified fractional sub-equation 
method and we present new Wick-type exact traveling wave solutions 
of Wick-type stochastic fractional RLW-Burgers equations. Finally, 
some conclusions are presented at the end of the paper in Section~\ref{sec:4}.


\section{Wick-type stochastic nonlinear Schr\"{o}dinger Equation}
\label{sec:2}

Let the space $(S)_{-1}^{n}$ consists of all formal expansions 
$F(z)=\sum_{\alpha} a_{\alpha} H_{\alpha}(z)$ 
with $a_{\alpha} \in \mathbb{R}^{n}$ and $z=(z_1, \ldots, z_n) 
\in \mathbb{C}^{n}$ (the set of all sequences of complex numbers). 


\subsection{Mathematical Computation scheme with Wick product}

We define the Wick product of two elements 
$F= \sum_{\alpha} a_{\alpha} H_{\alpha}$ and 
$G= \sum_{\beta} b_{\beta} H_{\beta} \in (S)_{-1}^n$ 
with $a_{\alpha}, b_{\beta} \in \mathbb{R}^n$ by
\begin{eqnarray*}
F \diamondsuit G = \sum_{\alpha, \beta} (a_{\alpha} , b_{\beta})
H_{\alpha+\beta}.
\end{eqnarray*}
For $F=\sum_{\alpha} a_{\alpha} H_{\alpha} \in (S)_{-1}^n$, with
$a_{\alpha} \in \mathbb{R}^n$, the Hermite transform of $F$,
denoted by  $\mathcal{H}(F)$ or $\tilde{F}(z)$ is defined by
\begin{eqnarray*}
\mathcal{H}(F)=\tilde{F}(z) = \sum a_{\alpha} z^{\alpha} \in
\mathbb{C}^{n}\ \   (\textrm{when} \mbox{~} \textrm{convergent}),
\end{eqnarray*}
where $z=(z_1, z_2, \ldots ) \in \mathbb{C}^{n}$ and $z^{\alpha} = (
z_1^{\alpha_1}, z_2^{\alpha_2} \ldots )$ for $\alpha=(\alpha_1,
\alpha_2 , \ldots ) \in J$, where $z_j^0 =1$.

By the above definition, for $F, G \in (S)_{-1}^{n}$,  we have
\begin{eqnarray*}
\widetilde{F \diamondsuit G}(z) = {\tilde F}(z) \dot {\tilde G}(z)
\end{eqnarray*}
for all $z$ such that ${\tilde F}(z) $ and ${\tilde G}(z)$ exist, where
the product on the right hand side of the above equation is the
complex bilinear product between two elements of
$\mathbb{C}^{n}$ defined by 
$(z_1^{\alpha_1}\ldots\, z_n^{\alpha_1}) \cdot
(z_1^{\alpha_2} \ldots\, z_n^{\alpha_2}) =\sum_{k=1}^{n} 
z_k^{\alpha_1} z_k^{\alpha_2}$, where $z_k^{\alpha_i} \in \mathbb{C}$.

We next present the main steps of finding exact solutions 
of the Wick-type stochastic nonlinear
partial differential equations.

Consider nonlinear stochastic
partial differential equation in the following form
\begin{eqnarray*}
\mathcal{P}(t,x,\partial_t, \nabla_x, U, z)=0 \label{eq1}
\end{eqnarray*}
where $\mathcal{P}$ is some given function, $U=U(t,x,z)$ is the
unknown solution, and the operators are 
$\partial_t = \frac{\partial}{\partial t}$, 
$\nabla_x =
\left(\frac{\partial}{\partial x_1}, \ldots,
\frac{\partial}{\partial x_d} \right)$ when $x=(x_1, \ldots, x_d)
\in \mathbb{R}^d $.\\

Suppose all functions and all products are interpreted 
as the Wick products and the Wick versions.
Eq. (\ref{eq1}) can be written in the form 
\begin{eqnarray}
\mathcal{P}^{\diamondsuit} (t,x, \partial_t, \nabla_x, U, \omega) =0.
\label{spde-1}
\end{eqnarray}
By taking the Hermite transform of Eq. (\ref{spde-1}),  
Eq. (\ref{spde-1}) convert Wick products into ordinary products
\begin{eqnarray}
{\tilde{ \mathcal{P}}} (t, x, \partial_t, \nabla_x, {\tilde U}, z)=0, 
\label{spde-2}
\end{eqnarray}
where ${\tilde U} =\mathcal{H}(U)$ is the Hermite transform of $U$
and $z=(z_1, z_2, \ldots) \in
\mathbb{C}^n$ \cite{benth,dai,chen}.\\

\noindent \textsf{Step 1.}  Combining the independent 
variables $x$ and $t$ into one variable $\eta$, we have
\begin{eqnarray}
{\tilde U}(x,t,z)= \textsf{u}(\eta), \eta =x-\int_0^t \omega(s,z)
ds. \label{tran1}
\end{eqnarray}
With the use of  the traveling wave variable (\ref{tran1}),
(\ref{spde-2}) is reduced to the ordinary differential equation for
$\textsf{u}=\textsf{u}(\eta)$
\begin{eqnarray}
\mathcal{Q}(\textsf{u}, \textsf{u}', \textsf{u}'',
\textsf{u}''', \ldots)=0, \label{ode}
\end{eqnarray}
where $\textsf{u}'=\frac{d \textsf{u}}{d \eta},
\textsf{u}''=\frac{d^2 \textsf{u}}{d \eta^2},
\textsf{u}'''=\frac{d^3 \textsf{u}}{d \eta^3}, \ldots$.

\noindent \textsf{Step 2.} Assume  that the solution of Eq.
(\ref{ode}) can be described by a polynomial in $\left\{
\textsf{F}(\eta)/\textsf{G}(\eta) \right\}$ as follows:
\begin{eqnarray}
\textsf{u}(\eta)= \sum_{i=0}^{m} A_i (t,z) \left\{
\frac{\textsf{F}(\eta)}{\textsf{G}(\eta)}\right\}^i,
\label{sol1}
\end{eqnarray}
where $\textsf{F}(\eta)$ and $\textsf{G}(\eta)$ 
satisfy the system with variable coefficients
\begin{eqnarray}
\cases{ \textsf{F}'(\eta) = p(t) \mbox{~} \textsf{F}(\eta), 
\textsf{G}'(\eta) = p(t) \mbox{~} \textsf{F}(\eta) + q(t) \mbox{~}
\textsf{G}(\eta),} \label{riccati-linear}
\end{eqnarray}
and $A_m(t,z), \ldots, A_1(t,z)$ and $A_0(t,z)$ are unknown 
coefficients to be computed later with $A_m(t,z) \neq 0$. 

The pole-order $m$ can be determined by considering the homogeneous 
balancing principle between the highest order derivatives of linear 
term and nonlinear term appearing in Eq. (\ref{ode}).
Now, the system (\ref{riccati-linear}) admits the following
ans$\ddot{{\rm a}}$tz:
\begin{eqnarray}
\left\{ \frac{\textsf{F}(\eta)}{\textsf{G}(\eta)} \right\} 
= \frac{p(t)-q(t)}{p(t)-q(t)\exp\{-(p(t)-q(t))\eta\}},
\label{new-sol1}
\end{eqnarray}
where $p(t)$ and $q(t)$ are arbitrary nonzero integrable 
functions with respect to parameter $t$ with $p(t) \neq q(t)$.

\noindent \textsf{Step 3.} Substituting the solution (\ref{sol1}) 
into Eq. (\ref{ode}) and collecting all terms with the same order 
of $\{\textsf{F}(\eta)/\textsf{G}(\eta)\}$ together, 
the left-hand side of Eq. (\ref{ode}) is converted into
another polynomial in $\{\textsf{F}(\eta) / \textsf{G}(\eta)\}$.
Equating each coefficient of this polynomial to zero, we get the
algebraic equations for coefficients $A_m (t,z), \ldots , A_1 (t,z),
A_0 (t,z)$ and $\omega(t,z)$. Further, the coefficients 
$A_m(t,z), \ldots, A_1(t,z), A_0(t,z)$, and $\omega(t,z)$ 
can be obtained by solving the algebraic equations, 
since the general solutions of the system (\ref{riccati-linear}) 
are known. Then by substituting $A_m(t,z), A_{m-1}(t,z),\ldots, 
A_1(t,z), A_0(t,z)$ and $\omega(t,z)$ into (\ref{sol1}), 
we can have more new exact solutions of Eq. (\ref{ode}).

\noindent \textsf{Step 4.} Some  explicit expressions for some undetermined
coefficients can derive from exact solutions of
(\ref{spde-2}) by substituting them into ${\tilde U}(t,x,z)$.
Moreover, under certain conditions, we take the inverse Hermite
transform $U=\mathcal{H}^{-1} ({\tilde U}) \in (S)_{-1}$ and thereby
obtain the Wick-type solutions $U$ of Eq. (\ref{spde-1}).

\begin{theorem}
\label{Th1}
Suppose $u(t,x,z)$ is a solution of (\ref{spde-2}) for $(t,x)$ 
in some bounded open set $\mathbf{G} \subset \mathbb{R} 
\times \mathbb{R}^d$, and for all $z \in \mathbb{K}_q(r)$, 
for some $q, r$. Also, assume that $u(t,x,z)$ 
and all its partial derivatives involved in
(\ref{spde-2}) are bounded for $(t,x,z) \in \mathbf{G} \times
\mathbb{K}_q(r)$, continuous with respect to $(t,x) \in \mathbf{G}$
for all $z \in \mathbb{K}_q(r)$ and analytic with respect to $z \in
\mathbb{K}_q(R)$, for all $(t,x) \in \mathbf{G}$. Then there exists
$U(t,x) \in H_{-1}$ such that $u(t,x,z) =  U(t,x)(z)$
for all $(t,x,z) \in \mathbf{G} \times \mathbb{K}_q(r)$ and $U(t,x)$
solves Eq. (\ref{spde-1}) in $H_{-1}$ in the strong sense.   
\end{theorem}


\subsection{Wick-type stochastic explicit solutions of equation (\ref{stochatic-schrodinger})}

By taking the Hermite transform in Eq. (\ref{stochatic-schrodinger}), we have
\begin{eqnarray}
i \tilde{\psi}_t(t,x,z) +\tilde{\alpha}(t,z) \tilde{
\psi}_{xx}(t,x,z) + \tilde{\beta}(t,z) \tilde{\psi}(t,x,z)
|\tilde{\psi}(t,x,z)|^2 + \tilde{\lambda}(t,z) \tilde{\psi}(t,x,z)
=0, \label{schrodinger-1}
\end{eqnarray}
where $z=(z_1,z_2, \ldots )\in \mathbb{C}^{n}$ is a
vector parameter. The condition $\alpha(t) \diamondsuit \beta(t) 
\diamondsuit \lambda(t) \neq 0$ yields $\tilde{\alpha}(t) 
\tilde{\beta}(t) \tilde{\lambda}(t) \neq 0$. To simplify 
Eq. (\ref{schrodinger-1}), we take
$\tilde\psi(t,x,z)=\psi(t,x,z)$, $\tilde\alpha(t,z)=\alpha(t,z)$,
$\tilde\beta(t,z)=\beta(t,z)$, and $\tilde\lambda(t,z)=\lambda(t,z)$.
Then, we have
\begin{eqnarray}
i \psi_t(t,x,z) +\alpha(t,z)\psi_{xx}(t,x,z) + \beta(t,z)
\psi(t,x,z) |\psi(t,x,z)|^2 + \lambda(t,z) \psi(t,x,z) =0.
\label{schrodinger-2}
\end{eqnarray}
In order to solve Eq. (\ref{schrodinger-2}), first we
substitute the transformation $\psi(t,x,z)=U(t,x,z) \textsf{e}^{i
\theta(t,z) }$ into Eq. (\ref{schrodinger-2}), and then Eq.
(\ref{schrodinger-2}) can be converted to the following equation
\begin{eqnarray}
i U_t(t,x,z) +\alpha(t,z)U_{xx}(t,x,z) + \beta(t,z) U^3(t,x,z) +
(\lambda(t,z)-\theta) U(t,x,z) =0, \label{schrodinger-3}
\end{eqnarray}
where $\theta$ is an integrable function of $t$. Now by using the
transformations
\begin{eqnarray}
U(t,x,z)= \textsf{u}(\eta), \eta= x-\int_0^t \omega(s,z) ds,
\label{transform-ode}
\end{eqnarray}
the Eq. (\ref{schrodinger-3}) can be rewritten in the following form
\begin{eqnarray}
-i \omega(t,z) \textsf{u}' +\alpha(t,z) \textsf{u}'' + \beta(t,z)
\textsf{u}^3 + (\lambda(t,z)-\theta(t,z)) \textsf{u} =0,
\label{schrodinger-ode}
\end{eqnarray}
where $\textsf{u}'= \frac{d \textsf{u}}{d \eta}$ and $\textsf{u}''=
\frac{d^2 \textsf{u}}{d \eta^2}$.
Now, let us find exact traveling wave solutions of 
Eq.~(\ref{stochatic-schrodinger}).
Suppose that the exact traveling wave solution 
of Eq. (\ref{schrodinger-ode}) can be
expressed in the form
\begin{eqnarray}
\textsf{u}(\eta) = A_1(t,z) \left\{
\frac{\textsf{F}(\eta)}{\textsf{G}(\eta)} \right\} +A_0(t,z) ,
\label{nls-solution}
\end{eqnarray}
where the ans$\ddot{{\rm a}}$tz $\{\textsf{F}(\eta)/\textsf{G}(\eta)\}$ is given by
\begin{eqnarray}
\left\{ \frac{\textsf{F}(\eta)}{\textsf{G}(\eta)} \right\} =
\frac{p(t)-q(t)}{p(t)-q(t)\exp\{-(p(t)-q(t))\eta\}},
\label{linear-sol}
\end{eqnarray}
where $\textsf{F}(\eta)$ and $\textsf{G}(\eta)$ satisfy system
(\ref{riccati-linear}), and $A_0(t,z), A_1(t,z)$ and $\omega(t,z)$
are time dependent functions to be determined later.
By taking the exact solution (\ref{nls-solution}) and the system
(\ref{riccati-linear}), we can obtain the derivatives $\textsf{u}'$
and $\textsf{u}''$ expressed via the ans$\ddot{{\rm a}}$tz
$\{\textsf{F}(\eta)/\textsf{G}(\eta)\}$. Substituting $\textsf{u},
\textsf{u}'$ and $\textsf{u}''$ into Eq. (\ref{schrodinger-ode}) and
equating to zero the expressions with the same degree of
$\{ \textsf{F}(\eta)/\textsf{G}(\eta) \}$ and solving the algebraic
equations with respect to the unknowns $A_1(t,z), A_0(t,z),
\lambda(t,z)$ and $\omega(t,z)$ by the help of Maple, we obtain
three nontrivial solution sets of coefficients as given below:
\begin{eqnarray}
\cases{ \theta(t,z)=\lambda(t,z)-2\alpha(t,z)(p(t)-q(t))^2,
\quad \omega(t,z)=-2 i \alpha(t,z)(p(t)-2q(t)), \cr A_0(t,z) =0, 
\quad 
A_1(t,z) = \pm p(t) \sqrt{-\frac{2\alpha(t,z)}{\beta(t,z)}}, } \label{coef-1}
\end{eqnarray}
\begin{eqnarray}
\cases{ \theta(t,z)=\lambda(t,z)-2\alpha(t,z)(p(t)-q(t))^2,
\quad \omega(t,z)=2 i \alpha(t,z)(2p(t)-q(t)), \cr A_0(t,z) =
\pm\frac{2\alpha(t,z)(p(t)-q(t))}{\beta(t,z)\sqrt{-\frac{2\alpha(t,z)}{\beta(t,z)}}},
\quad A_1(t,z) = \pm p(t) \sqrt{-\frac{2\alpha(t,z)}{\beta(t,z)}}, }
\label{coef-2}
\end{eqnarray}
\begin{eqnarray}
\cases{ \theta(t,z)=\lambda(t,z)- \frac{1}{2}
\alpha(t,z)(p(t)-q(t))^2, 
\quad \omega(t,z)= i \alpha(t,z)(p(t)+q(t)), \cr
A_0(t,z) =
\pm\frac{\alpha(t,z)(p(t)-q(t))}{\beta(t,z)\sqrt{-\frac{2\alpha(t,z)}{\beta(t,z)}}},
\quad A_1(t,z) = \pm p(t) \sqrt{-\frac{2\alpha(t,z)}{\beta(t,z)}}. }
\label{coef-3}
\end{eqnarray}
By substitutions (\ref{nls-solution}) and (\ref{linear-sol}) into the solution 
$\tilde{\psi}(t,x,z)=U(t,x,z)  \textsf{e}^{i\theta(t,z) }$, we
obtain new exact traveling wave solutions of Eq. (\ref{schrodinger-2}) 
according to each coefficient set (\ref{coef-1})--(\ref{coef-3}) as the
followings. Now, based on nontrivial coefficient set (\ref{coef-1}),  
exact traveling wave solution of Eq. (\ref{schrodinger-1}) can be written as
\begin{eqnarray}
\tilde{\psi}_1(t,x,z) =  \pm p(t)
\sqrt{-\frac{2\tilde{\alpha}(t,z)}{\tilde{\beta}(t,z)}} \left[
\frac{p(t)-q(t)}{p(t)-q(t)\exp\{-(p(t)-q(t))\eta(x,t,z)\}} \right]
\textsf{e}^{i \theta(t,z)}, \label{sol-1}
\end{eqnarray}
where $\eta(x,t,z)=x+2i \int_0^t \tilde{\alpha}(s,z)(p(s)-2q(s)) ds$
and
$\theta(t,z)=\tilde{\lambda}(t,z)-2\tilde{\alpha}(t,z)(p(t)-q(t))^2$.
Next, with the use of the nontrivial coefficient set (\ref{coef-2}), 
the exact traveling wave solution of Eq. (\ref{schrodinger-1}) is
\begin{eqnarray}
\tilde{\psi}_2(t,x,z) = \pm p(t)
\sqrt{-\frac{2\tilde{\alpha}(t,z)}{\tilde{\beta}(t,z)}} \left[
\frac{p(t)-q(t)}{p(t)-q(t)\exp\{-(p(t)-q(t))\eta(x,t,z)\}} \right]
\textsf{e}^{i \theta(t,z)} \pm
\frac{2\tilde{\alpha}(t,z)(p(t)-q(t))}{\tilde{\beta}(t,z)
\sqrt{-\frac{2\tilde{\alpha}(t,z)}{\tilde{\beta}(t,z)}}},
\label{sol-2}
\end{eqnarray}
where $\eta(x,t,z)=x - 2i \int_0^t \tilde{\alpha}(s,z)(2 p(s)- q(s))
ds$ and
$\theta(t,z)=\tilde{\lambda}(t,z)-2\tilde{\alpha}(t,z)(p(t)-q(t))^2$.
Finally, by using the nontrivial coefficient set (\ref{coef-3}), 
the exact traveling wave solution of Eq. (\ref{schrodinger-1}) is
\begin{eqnarray}
\tilde{\psi}_3(t,x,z) = \pm p(t)
\sqrt{-\frac{2\tilde{\alpha}(t,z)}{\tilde{\beta}(t,z)}} \left[
\frac{p(t)-q(t)}{p(t)-q(t)\exp\{-(p(t)-q(t))\eta(x,t,z)\}} \right]
\textsf{e}^{i \theta(t,z)} \pm
\frac{\tilde{\alpha}(t,z)(p(t)-q(t))}{\tilde{\beta}(t,z)
\sqrt{-\frac{2\tilde{\alpha}(t,z)}{\tilde{\beta}(t,z)}}},
\label{sol-3}
\end{eqnarray}
where $\eta(x,t,z)=x - i \int_0^t \tilde{\alpha}(s,z)(p(s) + q(s))
ds$ and $\theta(t,z)=\tilde{\lambda}(t,z)-\frac{1}{2}
\tilde{\alpha}(t,z)(p(t)-q(t))^2$.

In order to obtain white functional solutions for Eq.
(\ref{stochatic-schrodinger}), we use the Hermite transform and
Theorem~\ref{Th1}. Also, the property of generalized exponential functions
gives that there exists a bounded open set $\mathbf{G} \subset
\mathbb{R} \times \mathbb{R}^{+}$, $k<\infty$, $r>0$ such that the
exact solutions $\tilde{\psi}(t,x,z)$ of Eq. (\ref{schrodinger-1})
are uniformly bounded for $(t,x,z) \in \mathbf{G} \times
\mathbb{K}_q(r)$, continuous with respect to $(t,x) \in \mathbf{G}$
and analytic with respect to $z \in \mathbb{K}_k(r)$, for all $(x,t)
\in \mathbf{G}$. Also from Theorem~\ref{Th1}, there exist
$\Psi(t,x,z) \in (S)_{-1}$ such that $\tilde{\psi}(t,x,z) =
({\tilde\Psi }(t,x))(z)$ for all $(t,x,z) \in \mathbf{G} \times
\mathbb{K}_q(r)$ and $U(t,x)$ solves Eq.
(\ref{stochatic-schrodinger}) in $(S)_{-1}$. Hence, by applying the
inverse Hermite transform to solutions
(\ref{sol-1})--(\ref{sol-3}), we have the Wick-type solutions of Eq.
(\ref{stochatic-schrodinger}). Based on (\ref{sol-1}),
we get the white noise functional solution in the form
\begin{eqnarray}
\Psi_1 (t,x)=  \pm p(t) \diamondsuit
\sqrt{-\frac{2\alpha^*(t)}{\beta^*(t)}} \diamondsuit \left[
\frac{p(t)-q(t)}{p(t)-q(t)\diamondsuit
\exp^{\diamondsuit}\{-(p(t)-q(t))\diamondsuit \eta(x,t)\}} \right] 
\diamondsuit \textsf{e}^{\diamondsuit i \theta(t)}, \label{sol-1-1}
\end{eqnarray}
where $\eta(x,t)=x+2i \int_0^t \alpha^*(s) \diamondsuit (p(s)-2q(s)) ds$
and $\theta(t)=\lambda^*(t)-2 \alpha^*(t)\diamondsuit(p(t)-q(t))^2$.
Subsequently, from (\ref{sol-2}) and (\ref{sol-3}), we obtain the following
Wick-type solutions of Eq. (\ref{stochatic-schrodinger}):
\begin{eqnarray}
\Psi_2 (t,x)=  \pm p(t) \diamondsuit
\sqrt{-\frac{2\alpha^*(t)}{\beta^*(t)}} \diamondsuit \left[
\frac{p(t)-q(t)}{p(t)-q(t)\diamondsuit
\exp^{\diamondsuit}\{-(p(t)-q(t))\diamondsuit \eta(x,t)\}} \right] 
\diamondsuit \textsf{e}^{\diamondsuit i \theta(t)} \cr \pm \frac{2
\alpha^*(t)\diamondsuit(p(t)-q(t))}{\beta^*(t) \diamondsuit \sqrt{-\frac{2
\alpha^*(t)}{\beta^*(t)}}}, \label{sol-2-1}
\end{eqnarray}
where $\eta(x,t)=x-2i \int_0^t \alpha^*(s) \diamondsuit (2p(s)-q(s)) ds$
and $\theta(t)=\lambda^*(t)-2 \alpha^*(t)\diamondsuit(p(t)-q(t))^2$ and
\begin{eqnarray}
\Psi_3(x,t) = \pm p(t) \diamondsuit \sqrt{-\frac{2
\alpha^*(t)}{\beta^*(t)}} \diamondsuit \left[
\frac{p(t)-q(t)}{p(t)-q(t)\diamondsuit \exp^{\diamondsuit} \{-(p(t)-q(t))
\diamondsuit \eta(x,t)\}} \right] \diamondsuit \textsf{e}^{\diamondsuit 
i \theta(t)} \cr \pm \frac{\alpha^*(t)\diamondsuit(p(t)-q(t))}{\beta^*(t) 
\diamondsuit \sqrt{-\frac{2 \alpha^*(t)}{\beta^*(t)}}}, \label{sol-3-1}
\end{eqnarray}
where $\eta(x,t)=x - i \int_0^t \alpha^*(s)\diamondsuit(p(s) + q(s)) ds$
and $\theta(t)=\lambda^*(t)-\frac{1}{2} \alpha^*(t) \diamondsuit
(p(t)-q(t))^2$.

It should be mentioned that for different forms of
$\alpha^*(t), \beta^*(t)$ and $\lambda^*(t)$, we can obtain
different exact traveling wave solutions of Eq. 
(\ref{stochatic-schrodinger}) from the
solutions (\ref{sol-1-1})--(\ref{sol-3-1}).

\begin{example}
Take  $\tilde{\beta}(t,z) = \tilde{\alpha}(t,z)/c$ for a
constant parameter $c$ , we obtain the following solution of
Eq. (\ref{schrodinger-1}) based on  (\ref{sol-1}),
\begin{eqnarray}
\tilde{\psi}_{12}(t,x,z) =  \pm p(t) \sqrt{-2c} \left[
\frac{p(t)-q(t)}{p(t)-q(t)\exp\{-(p(t)-q(t))\eta(t,x,z))\}}\right]
\textsf{e}^{i \theta(t,z)}, \label{sol-12}
\end{eqnarray}
where $\eta(x,t,z)=x+2i \int_0^t \tilde{\alpha}(s,z)(p(s)-2q(s)) ds$
and
$\theta(t,z)=\tilde{\lambda}(t,z)-2\tilde{\alpha}(t,z)(p(t)-q(t))^2$. 
Further, we have the following solutions of Eq. (\ref{schrodinger-1}) 
from (\ref{sol-2}) and (\ref{sol-3}):
\begin{eqnarray}
\tilde{\psi}_{22}(t,x,z) =  \pm p(t) \sqrt{-2c} \left[
\frac{p(t)-q(t)}{p(t)-q(t)\exp\{-(p(t)-q(t))\eta(t,x,z))\}}\right]
\textsf{e}^{i \theta(t,z)} \pm \frac{2c(p(t)-q(t))}{
\sqrt{-2c}}, \label{sol-22}
\end{eqnarray}
where $\eta(x,t,z)=x-2 i \int_0^t \tilde{\alpha}(s,z)(2p(s)-q(s))
ds$ and
$\theta(t,z)=\tilde{\lambda}(t,z)-2\tilde{\alpha}(t,z)(p(t)-q(t))^2$,
\begin{eqnarray}
\tilde{\psi}_{32}(t,x,z) =  \pm p(t) \sqrt{-2c} \left[
\frac{p(t)-q(t)}{p(t)-q(t)\exp\{-(p(t)-q(t))\eta(t,x,z))\}}\right]
\textsf{e}^{i \theta(t,z)} \pm \frac{c(p(t)-q(t))}{ \sqrt{-2c}},
\label{sol-32}
\end{eqnarray}
where $\eta(x,t,z)=x- i \int_0^t \tilde{\alpha}(s,z)(p(s)+q(s)) ds$
and
$\theta(t,z)=\tilde{\lambda}(t,z)-\frac{1}{2}\tilde{\alpha}(t,z)(p(t)-q(t))^2$.
\end{example}

\begin{example}
It is noted that $\Psi(t,x)$ is the inverse Hermit
transformation of $\psi(t,x,z)$. If $\alpha^*(t) \ne 0$, then we can
obtain the following Wick-type solution 
of Eq. (\ref{stochatic-schrodinger}) from (\ref{sol-12}):
\begin{eqnarray}
\Psi_{13}(t,x) =  \pm  p(t) \sqrt{-2c} \diamondsuit \left[
\frac{p(t)-q(t)}{p(t)-q(t) \diamondsuit \exp^{\diamondsuit}\{-(p(t)-q(t))
\diamondsuit \eta(t,x))\}}\right] \diamondsuit \textsf{e}^{\diamondsuit i \theta(t)},
\label{sol-13}
\end{eqnarray}
where $\eta(x,t)=x+2i \int_0^t \alpha^*(s)\diamondsuit(p(s)-2q(s)) ds$
and $\theta(t)=\lambda^*(t)-2\alpha^*(t)\diamondsuit(p(t)-q(t))^2$.
Further, we obtain the following solutions from (\ref{sol-22}) and
(\ref{sol-32}):
\begin{eqnarray}
\Psi_{23}(t,x) =   \pm  p(t) \sqrt{-2c} \diamondsuit \left[
\frac{p(t)-q(t)}{p(t)-q(t) \diamondsuit \exp^{\diamondsuit}\{-(p(t)-q(t))
\diamondsuit \eta(t,x))\}}\right] \diamondsuit \textsf{e}^{\diamondsuit i \theta(t)} 
\pm \frac{2c(p(t)-q(t))}{ \sqrt{-2c}}, \label{sol-23}
\end{eqnarray}
where $\eta(x,t)=x-2 i \int_0^t \alpha^*(s)\diamondsuit(2p(s)-q(s)) ds$
and $\theta(t)=\lambda^*(t)-2 \alpha^*(t)\diamondsuit(p(t)-q(t))^2$,
\begin{eqnarray}
\Psi_{33}(t,x) =  \pm  p(t) \sqrt{-2c} \diamondsuit \left[
\frac{p(t)-q(t)}{p(t)-q(t) \diamondsuit \exp^{\diamondsuit}\{-(p(t)-q(t))
\diamondsuit \eta(t,x))\}}\right] \diamondsuit \textsf{e}^{\diamondsuit 
i \theta(t)} \pm \frac{c(p(t)-q(t))}{ \sqrt{-2c}}, \label{sol-33}
\end{eqnarray}
where $\eta(x,t)=x- i \int_0^t \alpha^*(s)\diamondsuit(p(s)+q(s)) ds$
and $\theta(t)=\lambda^*(t)-\frac{1}{2}
\alpha^*(t)\diamondsuit(p(t)-q(t))^2$.
\end{example}

\begin{example}
Let $f(t)$ be an  integrable function on
$\mathbb{R}_{+}$ and take $\alpha^*(t)=\alpha(t) + W(t)$ and
$\lambda^*(t)=\lambda(t) + W(t)$, where $W(t)$ denotes Gaussian white
noise, that is, $W(t) = \dot B(t) = dB(t)/dt$, 
$\{B(t) | t \in \mathbb{R}\}$ is a Brownian motion.
Further, taking the Hermite transform in $\alpha^*(t)$ and
$\lambda^*(t)$, we have $\alpha^*(t,z)=\alpha(t)+ \tilde{W} (t,z)$
and $\lambda^*(t)=\lambda(t)+ \tilde{W} (t,z)$, where $\tilde W
(t,z) = \sum_{k=1}^{\infty} \int_0^t \eta_k(s) ds z_k$ and $z=(z_1 ,
z_2 , \ldots ) \in \mathbb{C}^{n} $ is a parameter vector
and $\eta_k(s)$ is defined as in \cite{ycxie}. Now, based on  
(\ref{sol-13}), the stochastic soliton-like solution 
of Eq. (\ref{schrodinger-1}) can be obtained as
\begin{eqnarray}
\Psi_{14}(x,t,z) =  \pm  p(t) \sqrt{-2c} \left[
\frac{p(t)-q(t)}{p(t)-q(t)\exp\{-(p(t)-q(t)) \eta(t,x,z))\}}\right]
\textsf{e}^{ i \theta(t,z)}, \label{sol-14}
\end{eqnarray}
where $\eta(x,t,z)=x+2i \int_0^t
(\alpha(s)+\tilde{W}(s,z))(p(s)-2q(s)) ds$ and
\begin{equation}
\label{eq:theta}
\theta(t,z)=(\lambda(t)+\tilde{W}(t,z))-2
(\alpha(t)+\tilde{W}(t,z))(p(t)-q(t))^2. 
\end{equation}
Similarly, from (\ref{sol-23}) and (\ref{sol-33}) the solutions 
of Eq. (\ref{schrodinger-1}) can be obtained by
\begin{eqnarray}
\Psi_{24}(x,t,z) = \pm  p(t) \sqrt{-2c} \left[
\frac{p(t)-q(t)}{p(t)-q(t)\exp\{-(p(t)-q(t)) \eta(t,x,z))\}}\right]
\textsf{e}^{ i \theta(t,z)} \pm \frac{2c(p(t)-q(t))}{
\sqrt{-2c}}, \label{sol-24}
\end{eqnarray}
where $\eta(x,t,z)=x-2 i \int_0^t
(\alpha(s)+\tilde{W}(s,z))(2p(s)-q(s)) ds$, 
$\theta(t,z)$ is given by (\ref{eq:theta}) and
\begin{eqnarray}
\Psi_{34}(x,t,z) = \pm  p(t) \sqrt{-2c} \left[
\frac{p(t)-q(t)}{p(t)-q(t)\exp\{-(p(t)-q(t)) \eta(t,x,z))\}}\right]
\textsf{e}^{ i \theta(t,z)} \pm \frac{c(p(t)-q(t))}{
\sqrt{-2c}}, \label{sol-34}
\end{eqnarray}
where $\eta(x,t,z)=x- i \int_0^t
(\alpha(s)+\tilde{W}(s,z))(p(s)+q(s)) ds$ and
$$
\theta(t,z)=(\lambda(t)+\tilde{W}(t,z))-\frac{1}{2}
(\alpha(t)+\tilde{W}(t,z))(p(t)-q(t))^2.
$$
By non-random terms and random terms, the non-Wick-type solution
from (\ref{sol-24}) can be represented by
\begin{eqnarray}
\Psi_{15}(x,t) =  \pm  p(t) \sqrt{-2c} \left[
\frac{p(t)-q(t)}{p(t)-q(t)\exp\{-(p(t)-q(t)) \eta(t,x))\}}\right]
\textsf{e}^{ i \theta(t)}, \label{sol-15}
\end{eqnarray}
where $\eta(x,t)=x+2i \int_0^t (\alpha(s)+ \dot{B}(s))(p(s)-2q(s))
ds$ and $\theta(t)=(\lambda(t)+\dot{B}(t))-2
(\alpha(t)+\dot{B}(t))(p(t)-q(t))^2$. Similarly, from (\ref{sol-24}) and
(\ref{sol-34}), the solutions of Eq. (\ref{schrodinger-1}) are given by
\begin{eqnarray}
\Psi_{25}(x,t) = \pm  p(t) \sqrt{-2c} \left[
\frac{p(t)-q(t)}{p(t)-q(t)\exp\{-(p(t)-q(t)) \eta(t,x))\}}\right]
\textsf{e}^{ i \theta(t)} \pm \frac{2c(p(t)-q(t))}{ \sqrt{-2c}},
\label{sol-25}
\end{eqnarray}
where $\eta(x,t)=x-2 i \int_0^t (\alpha(s)+\dot{B}(s))(2p(s)-q(s))
ds$ and $\theta(t)=(\lambda(t)+\dot{B}(t))-2
(\alpha(t)+\dot{B}(t))(p(t)-q(t))^2$,
\begin{eqnarray}
\Psi_{35}(x,t,z) = \pm  p(t) \sqrt{-2c} \left[
\frac{p(t)-q(t)}{p(t)-q(t)\exp\{-(p(t)-q(t)) \eta(t,x))\}}\right]
\textsf{e}^{ i \theta(t)} \pm \frac{c(p(t)-q(t))}{ \sqrt{-2c}},
\label{sol-35}
\end{eqnarray}
where $\eta(x,t)=x- i \int_0^t (\alpha(s)+\dot{B}(s))(p(s)+q(s)) ds$
and $\theta(t,z)=(\lambda(t)+\dot{B}(t))-\frac{1}{2}
(\alpha(t)+\dot{B}(t))(p(t)-q(t))^2$.
\end{example}

\begin{remark}
By choosing appropriate values for physical parameters, 
the behaviors of the solitons-like traveling wave solution
(\ref{sol-25}) and (\ref{sol-35}) are shown graphically 
in Fig.~\ref{fig1} and Fig.~\ref{fig2}. 
\begin{figure}[ht!]
\begin{center}
(a) \epsfig{file=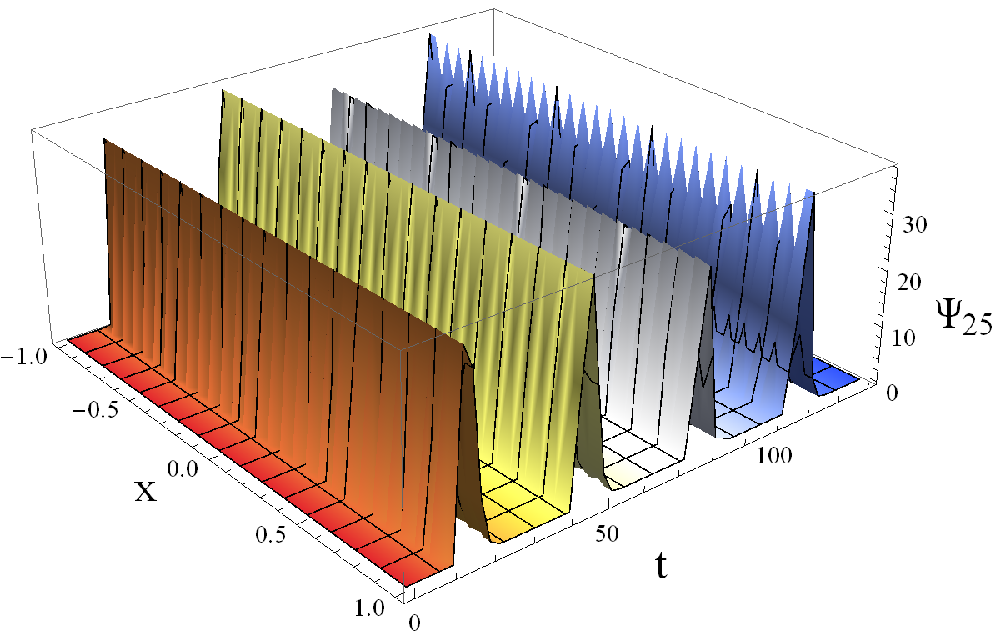, scale=0.5} \mbox{~~~}
(b) \epsfig{file=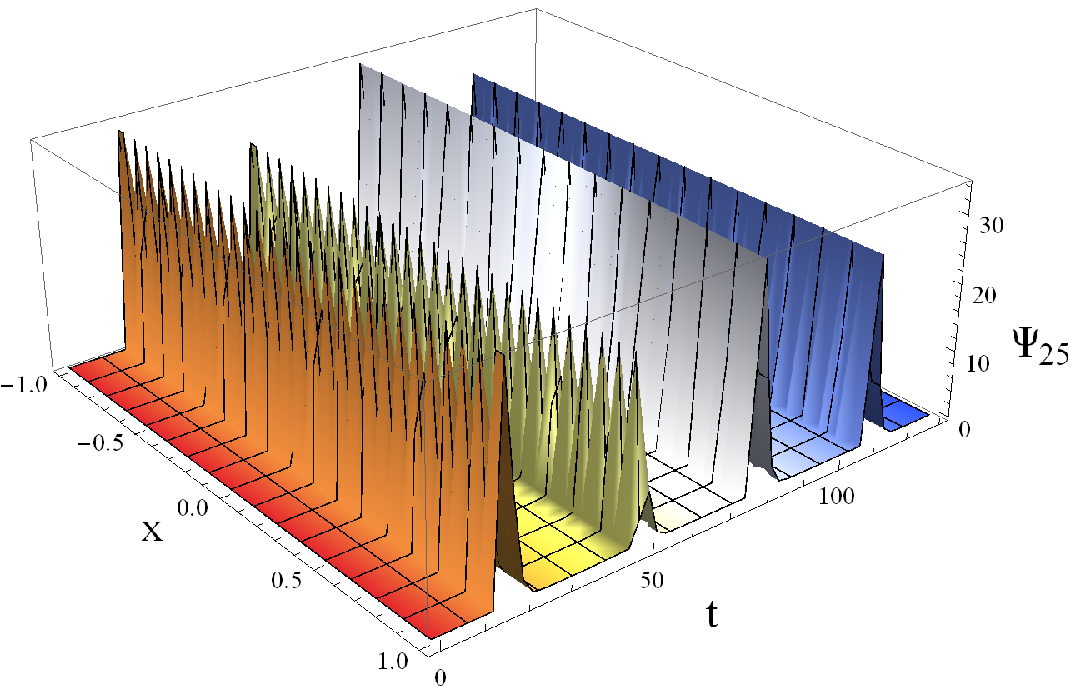, scale=0.5}
\end{center}
\vspace{-0.5cm}
\caption{Profiles of the non-Wick-type solution (\ref{sol-25}):
(a) $\dot{B}_t = i \cos(0.1t)$ and (b) $\dot{B}_t = 0$, 
under $p(t)= -1$, $q(t) = 0.2$, $c = -0.01$, 
$\alpha(t) = -1.5 i \cos(0.2 t)$, $\lambda(t) = i \cos(0.2 t)$.}  
\label{fig1}
\end{figure}
\begin{figure}[ht!]
\begin{center}
(a) \epsfig{file=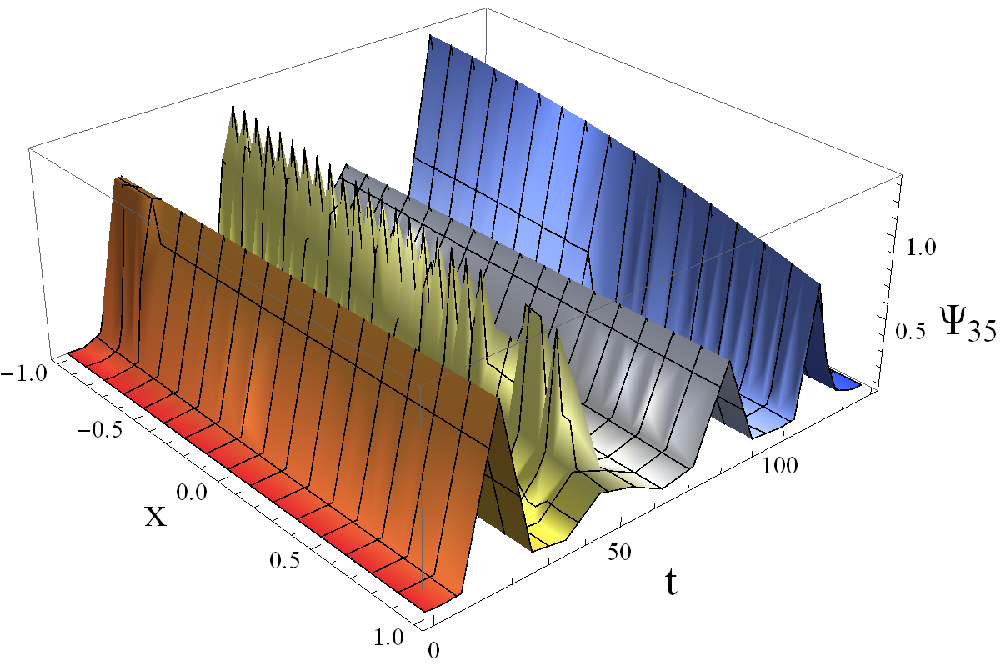, scale=0.55} \mbox{~~~}
(b) \epsfig{file=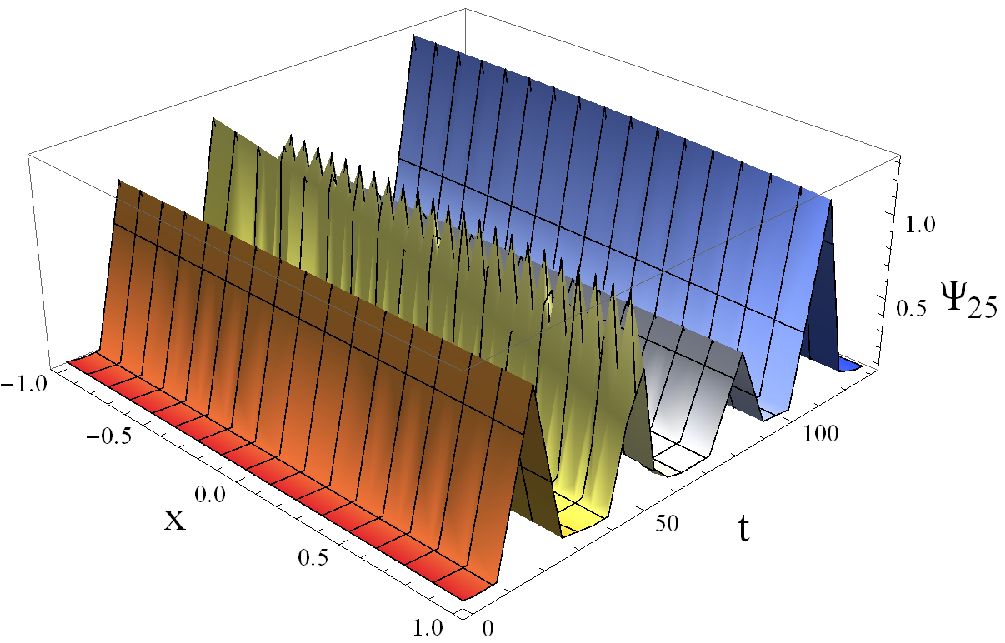, scale=0.55}
\end{center}
\vspace{-0.5cm}
\caption{Profiles of the non-Wick-type solution (\ref{sol-35}):
(a) $\dot{B}(t) = 0.25 i (\cos(0.1 t)+1.5 \sin(0.5 t))$; 
(b) $\dot{B}(t) = 0$, under $p(t)= -1$, $q(t) = 0.2$, $c = -0.01$, 
$\alpha(t) = -1.5 i \cos(0.2 t)$, $\lambda(t) = i \cos(0.2 t)$.} 
\label{fig2}
\end{figure}
In particular, Fig.~\ref{fig1} and \ref{fig2} are plotted 
for time-dependent coefficients $\alpha(t), \lambda(t)$,
without white noise effect and with white noise effect, respectively.
In Fig.~\ref{fig1}, with white noise functional, the white noise effect 
is generated in the part (a) as time changes, and without white noise functional, 
the noise is disappearing as time increases in the part (b). Fig.~\ref{fig2} 
performs the behaviors of the solitons-like traveling wave solution 
(\ref{sol-35}) or the white noise functionals.
\end{remark}


\section{Wick-type stochastic fractional RLW-Burgers equation}
\label{sec:3}

The Caputo fractional derivative of order $0<\alpha < 1$ is defined as follows \cite{kilbas}:
\begin{eqnarray}
{}^c D_t^{\alpha} f(t) = \frac{1}{\Gamma(1-\alpha)}  
\frac{\partial}{\partial \tau}
\int_0^t \frac{f(\tau)}{(t-\tau)^{\alpha}} d\tau .  
\label{frac-deri-1}
\end{eqnarray}
Some useful properties for the Caputo regularized partial derivative are:
\begin{eqnarray}
{}^c D_t^{\alpha} t^r
= \frac{\Gamma(1+r)}{\Gamma(1+r-\alpha)} t^{r-\alpha}, 
\label{frac-deri-2}
\end{eqnarray}
\begin{eqnarray}
{}^c D_t^{\alpha}(f(t)g(t)) 
= g(t) {}^c D_t^{\alpha} f(t) + f(t) {}^c D_t^{\alpha} g(t), 
\label{frac-deri-3}
\end{eqnarray}
\begin{eqnarray}
{}^c D_t^{\alpha} f[g(t)] = f'[g(t)] {}^c D_t^{\alpha} g(t).
\label{frac-deri-4}
\end{eqnarray}
In this section, we present the main steps of the mathematical
computation with the Caputo  fractional derivative for
obtaining exact solutions of fractional NPDEs. 


\subsection{Mathematical computation scheme with fractional order}

Consider the fractional improved Riccati equations with parameter
functions as follows:
\begin{eqnarray}
\label{riccati-frac}
\cases{ 
D_{\xi}^{\alpha} \Psi(\xi) = \lambda(t)  \Psi(\xi),
\cr D_{\xi}^{\alpha} \Phi(\xi) = \lambda(t) \Psi(\xi)
+\mu(t) \Phi(\xi), }
\end{eqnarray}
where $\lambda(t)$ and $\mu(t)$ denote arbitrary integrable functions,
$D_{\xi}^{\alpha} \Psi(\xi)$ and $D_{\xi}^{\alpha} \Phi(\xi)$ denote
the Caputo fractional derivative of fractional-order $\alpha$ for 
$\Psi(\xi)$ and $\Phi(\xi)$, and ${}^c  D_{\xi}^{\alpha}=D_{\xi}^{\alpha}$ 
for using a simple notation. In order to obtain the
solutions for Eq. (\ref{riccati-frac}), we consider
$\Psi(\xi)=\psi(\eta)$ and $\Phi(\xi)=\phi(\eta)$ with a nonlinear
fractional complex transformation $\eta= \xi^{\alpha}
/\Gamma(1+\alpha)$. By using Eq. (\ref{frac-deri-2}) and Eq.
(\ref{frac-deri-4}), that is, $D_{\xi}^{\alpha} \Psi(\xi) =
D_{\xi}^{\alpha} \psi(\eta) = \psi'(\eta)D_{\xi}^{\alpha} \eta=
\psi'(\eta)$ and $D_{\xi}^{\alpha} \Phi(\xi) = D_{\xi}^{\alpha}
\phi(\eta) = \phi'(z)D_{\xi}^{\alpha} \eta = \phi'(\eta)$, the
fractional improved Riccati equations (\ref{riccati-frac}) can be
transformed to the improved Riccati equation in the following form:
\begin{eqnarray}
\label{riccati-ode}
\cases{  \psi'(\eta) = \lambda(t) \mbox{~} \psi(\eta), \cr
\phi'(\eta) = \lambda(t) \mbox{~} \psi(\eta)+ \mu(t) \phi(\eta).}
\end{eqnarray}

From the general solutions of Eq. (\ref{riccati-ode}), we construct
the ans$\ddot{{\rm a}}$tz in the following form:
\begin{eqnarray}
\left\{\frac{\psi(\eta)}{\phi(\eta)}\right\} =
\frac{\lambda(t)-\mu(t)}{\lambda(t)-
\mu(t)\exp\{-(\lambda(t)-\mu(t)) \eta \}}, 
\label{riccati-ode-sol}
\end{eqnarray}
and by substituting the nonlinear fractional transformation
$\eta= \xi^{\alpha} /\Gamma(1+\alpha)$, the ans$\ddot{{\rm a}}$tz
(\ref{riccati-ode-sol}) can be rewritten in the following form:
\begin{eqnarray}
\left\{\frac{\Psi(\xi)}{\Phi(\xi)}\right\} =
\frac{\lambda(t)-\mu(t)}{\lambda(t)-
\mu(t)\exp\{-(\lambda(t)-\mu(t))\frac{\xi^{\alpha}}{\Gamma(1+\alpha)}
\}}. \label{riccati-frac-sol}
\end{eqnarray}
Suppose that the fractional Wick-type stochastic NPDE with the
independent variables $t, x_1, x_2, \ldots, x_n$, is given by
\begin{eqnarray}
\mathcal{P}^{\diamondsuit}\left(u, D_t^{\diamondsuit \alpha} u,
D_{x_1}^{\diamondsuit\alpha} u, D_{x_2}^{\diamondsuit\alpha} u, 
\ldots,
D_t^{\diamondsuit 2\alpha} u, D_{x_1}^{\diamondsuit 2\alpha} u,
D_{x_2}^{\diamondsuit 2\alpha} u, \ldots \right)=0, \label{wick-pde1}
\end{eqnarray}
where $D_t^{\diamondsuit \alpha} u, D_{x_1}^{\diamondsuit\alpha} u,
D_{x_2}^{\diamondsuit\alpha} u, \ldots, D_t^{\diamondsuit 2\alpha} u,
D_{x_1}^{\diamondsuit 2\alpha} u, D_{x_2}^{\diamondsuit 2\alpha} u, \ldots$,
are the Caputo fractional derivatives of $u$ in the
Wick-type sense with respect to $t, x_1, x_2, \ldots, x_n$ and
$u=u(t, x_1 , x_2, \ldots, x_n)$. $\mathcal{P}$ is a polynomial in
$u$ and its various partial fractional derivatives.\\

\noindent \textsf{Step 1.}  Suppose that traveling wave
transformation is given by
\begin{eqnarray}
\cases{ u(t, x_1, x_2, \ldots, x_n)=U(\zeta), \cr \zeta=
\int_0^{T} c(s) ds +\frac{k_1
x_1^{\alpha}}{\Gamma(1+\alpha)}+\frac{k_2
x_2^{\alpha}}{\Gamma(1+\alpha)} + \cdots + \frac{k_n
x_n^{\alpha}}{\Gamma(1+\alpha)},  T
= \frac{ t^{\alpha}}{\Gamma(1+\alpha)} } \label{twv}
\end{eqnarray}
where $k_1, k_2, \ldots, k_n$ represent nonzero constants 
and $c(t)$ is an integrable function of $t$. Then, 
by the property (\ref{frac-deri-4}) 
of the Caputo fractional derivative and by taking
$u=\tilde{U}(\zeta)$, $D_t^{\alpha} u =c(T) \tilde{U}_{\zeta}$,
$D_{x_1}^{\alpha} u =k_1 \tilde{U}_{\zeta}$, $D_{x_2}^{\alpha} u =k_2
\tilde{U}_{\zeta}$, $\ldots$, $D_t^{2\alpha} u =c^2(T)
\tilde{U}_{\zeta\zeta}$, $D_{x_1}^{2\alpha} u =k_1^2
\tilde{U}_{\zeta\zeta}$, $D_{x_2}^{2\alpha} u = k_2^2
\tilde{U}_{\zeta\zeta}$, $\ldots$,  the  Eq. (\ref{wick-pde1}) can be
reduced into the following fractional ordinary differential equation
with respect to traveling wave variable $\zeta$:
\begin{eqnarray}
\tilde{\mathcal{P}}(\tilde{U}, c \tilde{U}', k_1 \tilde{U}', k_2
\tilde{U}' \ldots, c^2 \tilde{U}'', k_1^2 \tilde{U}'', k_2^2
\tilde{U}'' \ldots, z)=0, \label{frac-pde2}
\end{eqnarray}
where $c=c(T),\tilde{U}'=\tilde{U}_{\zeta},
\tilde{U}''=\tilde{U}_{\zeta\zeta}, \ldots$, and
$\tilde{U}=\mathcal{H}(U)$ is the Hermite transform of $U$ and
$z=(z_1, z_2 , \ldots)$ is a vector of all sequences of complex
numbers. Suppose that we can find the solution $u=\tilde{U}(t,x,z)$
of Eq. (\ref{frac-pde2}) for $z \in \mathbb{K}_m(n)$, where
$\mathbb{K}_m(n)=\{(z_1, z_2, \ldots ) \in \mathbb{C}^{n}
\mbox{~}{\rm and}\mbox{~} \sum_{\alpha} |z^{\alpha}|^2
(2 N)^{m \alpha} <n^2\}$ for some integers $m, n$, $x=(x_1,
x_2, \ldots, x_n)$ and $z=(z_1, z_2, \ldots, z_n)$.\\

\noindent \textsf{Step 2.} Suppose that the solution of Eq.
(\ref{frac-pde2}) can be expressed by a polynomial in
$\left\{\psi(\zeta)/\phi(\zeta)\right\}$ as given below:
\begin{eqnarray}
\tilde{U}(\zeta)= \sum_{i=0}^m A_i(t) \left\{
\frac{\psi(\zeta)}{\phi(\zeta)} \right\}^i, \label{solution}
\end{eqnarray}
where $\left\{\psi(\zeta)/\phi(\zeta)\right\}$ is the ans$\ddot{{\rm a}}$tz 
(\ref{riccati-ode-sol}) and $\{A_{i}(t) \}_{ i=0}^{m}$ are
unknown coefficients to be computed later, $A_{m}(t) \neq 0$. The
pole-order $m$ can be determined by considering the homogeneous
balancing principle between the highest order linear term and the
highest order nonlinear term that occurs in Eq. (\ref{frac-pde2}).\\

\noindent \textsf{Step 3.} Substituting the solution
(\ref{solution}) into Eq. (\ref{frac-pde2}), collecting all terms
with the same order of $\left\{\psi(\zeta)/\phi(\zeta)\right\}$
together, the left-hand side of Eq. (\ref{frac-pde2}) changed  into
another polynomial in $\left\{\psi(\zeta)/\phi(\zeta)\right\}$.
Equating each coefficient of this resulting polynomial to zero, we
get a set of  algebraic equations for coefficients $\left\{A_i(t)
\right\}_{i=0}^m, \left\{k_i \right\}_{i=1}^n$ and $c(t)$.\\

\noindent \textsf{Step 4.} Solving the resulting set of obtained 
algebraic equations in Step 3 and combining the solution
(\ref{solution}), the traveling wave transformation (\ref{twv}) and
the fractional ans$\ddot{\rm a}$tz (\ref{riccati-frac-sol}), under
certain conditions, we can take the inverse Hermite transform
$\tilde{U} =\mathcal{H}^{-1}(u) \in (\mathcal{S})_{-1}$ and hence 
we obtain the Wick-type solution $u$ of the original fractional Wick-type
stochastic NPDE (\ref{wick-pde1}) \cite{holden}.


\subsection{Wick-type fractional exact traveling wave solutions of equation (\ref{rlw-wick})}

Now, by applying the Hermite transform, the fractional
Wick-type stochastic RLW-Burgers equation is converted into a
deterministic fractional NPDE, which is expressed in the following form:
\begin{eqnarray}
\frac{\partial^{\alpha} \tilde{V}}{\partial t^{\alpha}} +
\tilde{P}(t,z)\frac{\partial^{\alpha} \tilde{V}}{\partial
x^{\alpha}} +\tilde{Q}(t,z)\tilde{V} \frac{\partial^{\alpha}
\tilde{V}}{\partial x^{\alpha}}
+\tilde{R}(t,z)\frac{\partial^{2\alpha} \tilde{V}}{\partial
x^{2\alpha}}+\tilde{S}(t,z)\frac{\partial^{3\alpha}
\tilde{V}}{\partial t^{\alpha}
\partial x^{2\alpha}} =0, 
\label{rlw-hermite}
\end{eqnarray}
where $z=(z_1, z_2, \ldots) \in \mathbb{C}^{n}$ is a vector
parameter. To derive exact traveling wave solutions of Eq.
(\ref{rlw-hermite}), we consider the solution form
$\tilde{V}=\tilde{V}(x,t,z)=\textsf{y}(\zeta(x,t,z))$ with the
traveling wave variable
\begin{eqnarray}
\zeta(x,t,z)= k X + \int_0^{T} c(s,z) ds, \label{rlw-twv}
\end{eqnarray}
where $X=\frac{x^{\alpha}}{\Gamma(1+\alpha)},
T=\frac{t^{\alpha}}{\Gamma(1+\alpha)}$, $k$ is a nonzero constant
and $c(s,z)$ is a nonzero function of the indicated variables to be
calculated  later. By using the traveling wave variable
(\ref{rlw-twv}), Eq. (\ref{rlw-hermite}) can be converted into the form
\begin{eqnarray}
(c(T,z)+k p(t,z))\textsf{y}' + k q(t,z) \textsf{y}\textsf{y}'
+k^2r(t,z) \textsf{y}'' +kc(T,z)s(t,z)\textsf{y}''' =0 ,
\label{rlw-hermite-1}
\end{eqnarray}
where $\textsf{y}'= \frac{d \tilde{V}}{d \zeta}, \textsf{y}''=
\frac{d^2 \tilde{V}}{d \zeta^2}, \textsf{y}'''= \frac{d^3
\tilde{V}}{d \zeta^3}$ and $p(t,z)=\tilde{P}(t,z),
q(t,z)=\tilde{Q}(t,z), r(t,z)=\tilde{R}(t,z),
s(t,z)=\tilde{S}(t,z).$ Integrating Eq. (\ref{rlw-hermite-1}) with
respect to $\zeta$ once, we can get
\begin{eqnarray}
(c(T,z)+k p(t,z))\textsf{y} + \frac{1}{2} k q(t,z) \textsf{y}^2 +k^2
r(t,z) \textsf{y}' +k^2 c(T,z)s(t,z)\textsf{y}'' =0.
\label{rlw-hermite-2}
\end{eqnarray}
Here the arbitrary integral constant is assumed to be zero.

To determine the pole-order $m$ of the exact traveling wave solution 
of Eq.~(\ref{rlw-hermite-2}), by balancing the highest order linear 
term $\textsf{y}''$ and the highest order nonlinear term $\textsf{y}^2$
in Eq.~(\ref{rlw-hermite-2}), we obtain $N+2=2N$, which gives $N=2$.
Now, we assume that the second-order pole exact traveling wave solution of Eq.
(\ref{rlw-hermite-2}) can be expressed by the following form:
\begin{eqnarray}
\textsf{y}(\zeta)= \sum_{i=0}^2 A_i(t,z)
\left\{\frac{\psi(\zeta)}{\phi(\zeta)}\right\}^i.
\label{rlw-hermite-sol}
\end{eqnarray}
Substituting the solution
(\ref{rlw-hermite-sol}) into Eq. (\ref{rlw-hermite-2}), collecting
all the terms with the same power of $\{\psi(\zeta)/\phi(\zeta)\}$
together and equating each coefficient of this polynomial to zero,
we obtain a set of algebraic equations for $A_0(t,z), A_1(t,z),
A_2(t,z)$ and $c(t,z)$. Further, by solving the algebraic equations
and substituting the coefficients of nontrivial solutions and the 
traveling wave variable (\ref{rlw-twv}) into Eq. (\ref{rlw-hermite-sol}), 
we obtain six exact traveling wave solutions of Eq. (\ref{rlw-hermite}).
Let $\tau_{\alpha}=[\tau \Gamma(1+\alpha)]^{1/\alpha}$,
$B_{1}(t,z)=r(t,z)-s(t,z)\mu(t,z)$, 
$B_{2}(t,z)=r(t,z)-2s(t,z)\mu(t,z)$,
$C_{1}(t,z)=r(t,z)-6s(t,z)\mu(t,z)$,
$C_{2}(t,z)=r(t,z)-5s(t,z)\mu(t,z)$,
$D_{1}(t,z)=r(t,z)+4s(t,z)\mu(t,z)$,
$D_{2}(t,z)=r(t,z)-2s(t,z)\mu(t,z)$, 
and $k$ be a nonzero constant.

The first exact traveling wave solution with a relation
$\lambda(t,z)=(r(t,z)-s(t,z) \mu(t,z))/s(t,z)$ 
is expressed by the form
\begin{eqnarray}
\tilde{V}_1(x,t,z) = \frac{12k B_{1}(t,z)B_{2}^2(t,z)} {q(t,z)s(t,z)
\left[B_{1}(t,z)-s(t,z) \mu(t,z)\exp\{-\frac{B_2(t,z) }{s(t,z)}
\zeta_1(x,t,z) \} \right]}  \mbox{~~~~~} \cr -\frac{12k B_{1}^2(t,z)
B_{2}(t,z)} {q(t,z)s(t,z) \left[B_{1}(t,z)-s(t,z)
\mu(t,z)\exp\{-\frac{B_2(t,z) }{s(t,z)} \zeta_1(x,t,z) \} \right]^2}, 
\label{frac-sol-1}
\end{eqnarray}
where $\zeta_1(x,t,z) = \frac{k x^{\alpha}}{\Gamma(1+\alpha)} +
\int_0^{\frac{t^{\alpha}}{\Gamma(1+\alpha)}} c(\tau,z) d\tau$ and 
$$
c(\tau,z) = -k ( 4 k \mu^2(\tau_{\alpha},z) s(\tau_{\alpha},z) -4k
r(\tau_{\alpha},z) \mu(\tau_{\alpha},z) + p(\tau_{\alpha},z) + k
r^2(\tau_{\alpha},z)/s(\tau_{\alpha},z)).
$$

The second exact traveling wave solution with relation 
$\lambda(t,z)=-(r(t,z)-6 s(t,z) \mu(t,z))/s(t,z)$ is given by
\begin{eqnarray}
\tilde{V}_2(x,t,z) = -\frac{3k C_1^2(t,z)
C_2(t,z)}{4q(t,z)s(t,z)\left[ C_{1}(t,z)+s(t,z)
\mu(t,z)\exp\{-\frac{C_2(t,z) }{s(t,z)} \zeta_2(x,t,z) \} \right]^2}, 
\label{frac-sol-2}
\end{eqnarray}
where $\zeta_2(x,t,z) = \frac{k x^{\alpha}}{\Gamma(1+\alpha)} +
\int_0^{\frac{t^{\alpha}}{\Gamma(1+\alpha)}} c(\tau,z) d\tau$ and 
$$
c(\tau,z) =k( 12 k \mu^2(\tau_{\alpha},z) s(\tau_{\alpha},z) - 12 k
r(\tau_{\alpha},z)\mu(\tau_{\alpha},z) -8 p(\tau_{\alpha},z) +3 k
r^2(\tau_{\alpha},z)/ s(\tau_{\alpha},z) )/8.
$$

The third exact traveling wave solution with relation
$\lambda(t,z)=(r(t,z)+4s(t,z) \mu(t,z))/(6s(t,z))$ can be expressed by
\begin{eqnarray}
\tilde{V}_3(x,t,z) = \frac{2k D_1(t,z) D_2^2(t,z)}{q(t,z)s(t,z)
\left[ D_{1}(t,z)-6s(t,z) \mu(t,z)\exp\{-\frac{D_2(t,z)}{
6s(t,z)}\zeta_3(x,t,z) \} \right]} \mbox{~~~~~}\cr -\frac{k D_1^2(t,z)
D_2(t,z)}{3q(t,z)s(t,z) \left[ D_{1}(t,z)-6s(t,z)
\mu(t,z)\exp\{-\frac{D_2(t,z) }{6s(t,z)}\zeta_3(x,t,z) \} \right]^2}, 
\label{frac-sol-3}
\end{eqnarray}
where $\zeta_3(x,t,z) = \frac{k x^{\alpha}}{\Gamma(1+\alpha)} +
\int_0^{\frac{t^{\alpha}}{\Gamma(1+\alpha)}} c(\tau,z) d\tau$, $k$
an arbitrary constant, and
$$
c(\tau,z) = - k( 4 k \mu^2(\tau_{\alpha},z)
s(\tau_{\alpha},z)  - 4k r(\tau_{\alpha},z) \mu(\tau_{\alpha},z)+6
p(\tau_{\alpha},z)+ k r^2(\tau_{\alpha},z)/ s(\tau_{\alpha},z))/6.
$$

The fourth exact traveling wave solution with relation
$\lambda(t,z)=(r(t,z)-s(t,z) \mu(t,z))/s(t,z)$ is expressed by
\begin{eqnarray}
\tilde{V}_4(x,t,z) = -\frac{2k B_2^2(t,z)}{s(t,z)\mu(t,z)} \frac{12k
B_{1}(t,z)B_{2}^2(t,z)} {q(t,z)s(t,z) \left[B_{1}(t,z)-s(t,z)
\mu(t,z)\exp\{-\frac{B_2(t,z) }{s(t,z)} \zeta_4(x,t,z) \} \right]}
\mbox{~}\cr -\frac{12k B_{1}^2(t,z) B_{2}(t,z)} {q(t,z)s(t,z)
\left[B_{1}(t,z)-s(t,z) \mu(t,z)\exp\{-\frac{B_2(t,z) }{s(t,z)}
\zeta_4(x,t,z) \} \right]^2}, 
\label{frac-sol-4}
\end{eqnarray}
where $\zeta_4(x,t,z) = \frac{k x^{\alpha}}{\Gamma(1+\alpha)} +
\int_0^{\frac{t^{\alpha}}{\Gamma(1+\alpha)}} c(\tau,z) d\tau$ and 
$$
c(\tau,z) = k ( 4 k \mu^2(\tau_{\alpha},z) s(\tau_{\alpha},z) -4k
r(\tau_{\alpha},z) \mu(\tau_{\alpha},z) - p(\tau_{\alpha},z) + k
r^2(\tau_{\alpha},z)/s(\tau_{\alpha},z)).
$$

The fifth exact traveling wave solution with relation 
$\lambda(t,z)=-(r(t,z)-6 s(t,z) \mu(t,z))/s(t,z)$ is given by
\begin{eqnarray}
\tilde{V}_5(x,t,z) = \frac{3k B_2^2(t,z)}{4q(t,z)s(t,z)} -\frac{3k
C_1^2(t,z) C_2(t,z)}{4q(t,z)s(t,z)\left[ C_{1}(t,z)+s(t,z)
\mu(t,z)\exp\{-\frac{C_2(t,z) }{s(t,z)} \zeta_5(x,t,z)\} \right]^2}, 
\label{frac-sol-5}
\end{eqnarray}
where $\zeta_5(x,t,z) = \frac{k x^{\alpha}}{\Gamma(1+\alpha)} +
\int_0^{\frac{t^{\alpha}}{\Gamma(1+\alpha)}} c(\tau,z) d\tau$ and 
$$
c(\tau,z) = -k ( 12 k \mu^2(\tau_{\alpha},z) s(\tau_{\alpha},z) - 12
k r(\tau_{\alpha},z)\mu(\tau_{\alpha},z) +8 p(\tau_{\alpha},z) +3 k
r^2(\tau_{\alpha},z)/ s(\tau_{\alpha},z) )/8.
$$

Finally, the sixth exact traveling wave solution with relation
$\lambda(t,z)=(r(t,z)+4s(t,z) \mu(t,z))/(6s(t,z))$ is given by
\begin{eqnarray}
\tilde{V}_6(x,t,z) = \frac{2k D_1(t,z) D_2^2(t,z)}{q(t,z)s(t,z)
\left[ D_{1}(t,z)-6s(t,z) \mu(t,z)\exp\{-\frac{D_2(t,z)}{
6s(t,z)}\zeta_6(x,t,z) \} \right]}  \mbox{~~~~~}\cr -\frac{k D_1^2(t,z)
D_2(t,z)}{3q(t,z)s(t,z) \left[ D_{1}(t,z)-6s(t,z)
\mu(t,z)\exp\{-\frac{D_2(t,z) }{6s(t,z)}\zeta_6(x,t,z) \} \right]^2}, 
\label{frac-sol-6}
\end{eqnarray}
where $\zeta_6(x,t,z) = \frac{k x^{\alpha}}{\Gamma(1+\alpha)} +
\int_0^{\frac{t^{\alpha}}{\Gamma(1+\alpha)}} c(\tau,z) d\tau$ and
$$
c(\tau,z) =k ( 4 k \mu^2(\tau_{\alpha},z) s(\tau_{\alpha},z)  - 4k
r(\tau_{\alpha},z) \mu(\tau_{\alpha},z) -6 p(\tau_{\alpha},z)+ k
r^2(\tau_{\alpha},z)/ s(\tau_{\alpha},z))/6.
$$

In order to obtain the white noise functional solutions of Eq.
(\ref{rlw-wick}), we need to use the inverse Hermite transform and
Theorem 4.1.1 in \cite{holden}. Then we can obtain the Wick-type
exact traveling wave solutions as the white noise functional 
solutions of Eq. (\ref{rlw-wick}) based on the solutions
(\ref{frac-sol-1})--(\ref{frac-sol-6}) in the following form:
\begin{eqnarray}
V_1(x,t) = \frac{12k B_{1}(t)\diamondsuit B_{2}^{\diamondsuit 2}(t)}{
Q(t)\diamondsuit S(t) \diamondsuit \left[B_{1}(t)-S(t)\diamondsuit
\mu(t)\diamondsuit \exp^{\diamondsuit}\{-\frac{B_2(t)}{S(t)}\diamondsuit
\Xi_1(x,t) \}\right]}  \mbox{~}\cr -\frac{12k B_{1}^{\diamondsuit 2}(t) 
\diamondsuit B_{2}(t)}{Q(t)\diamondsuit S(t)\diamondsuit \left[B_{1}(t)-S(t)
\diamondsuit \mu(t)\diamondsuit \exp^{\diamondsuit }\{-\frac{B_2(t)}{S(t)}
\diamondsuit \Xi_1(x,t) \}\right]^{\diamondsuit 2}}, \label{wick-sol-1}
\end{eqnarray}
where $\Xi_1(x,t) = \frac{k x^{\alpha}}{\Gamma(1+\alpha)} +
\int_0^{\frac{t^{\alpha}}{\Gamma(1+\alpha)}} c(\tau) d\tau$ and
$$
c(\tau)=-k(4k \mu^{\diamondsuit 2}(\tau_{\alpha})\diamondsuit
S(\tau_{\alpha})-4k R(\tau_{\alpha})\diamondsuit \mu(\tau_{\alpha})
+P(\tau_{\alpha}) +k R^{\diamondsuit 2}(\tau_{\alpha})/S(\tau_{\alpha})), 
$$
\begin{eqnarray}
V_2(x,t) = -\frac{3k C_1^{\diamondsuit 2}(t) \diamondsuit C_2(t)}{4Q(t)
\diamondsuit S(t) \diamondsuit \left[ C_{1}(t)+S(t) \diamondsuit \mu(t)
\diamondsuit\exp^{\diamondsuit}\{-\frac{C_2(t) }{S(t)}\diamondsuit \Xi_2(x,t) 
\} \right]^{\diamondsuit 2}},
\label{wick-sol-2}
\end{eqnarray}
where $\Xi_2(x,t) = \frac{k x^{\alpha}}{\Gamma(1+\alpha)} +
\int_0^{\frac{t^{\alpha}}{\Gamma(1+\alpha)}} c(\tau) d\tau$ and
$$
c(\tau)= k(12k \mu^{\diamondsuit 2}(\tau_{\alpha})\diamondsuit
S(\tau_{\alpha})-12 k R(\tau_{\alpha})\diamondsuit \mu(\tau_{\alpha})
-8P(\tau_{\alpha})+3 k R^{\diamondsuit 2}(\tau_{\alpha})/
S(\tau_{\alpha}) )/8,
$$ 
\begin{eqnarray}
V_3(x,t) = \frac{2k D_1(t) \diamondsuit D_2^{\diamondsuit 2}(t)}{Q(t)
\diamondsuit S(t) \diamondsuit \left[ D_{1}(t)-6S(t) \diamondsuit
\mu(t) \diamondsuit \exp^{\diamondsuit}\{-\frac{D_2(t) }{6S(t)} 
\diamondsuit \Xi_3(x,t)\} \right] } \mbox{~~~~}\cr -\frac{k D_1^{\diamondsuit 2}(t) 
\diamondsuit D_2(t)}{3Q(t)\diamondsuit S(t) \diamondsuit \left[ D_{1}(t)-6S(t) 
\diamondsuit \mu(t) \diamondsuit \exp^{\diamondsuit}\{-\frac{D_2(t) }{6S(t)} 
\diamondsuit \Xi_3(x,t)\} \right]^{\diamondsuit 2}}, 
\label{wick-sol-3}
\end{eqnarray}
where $\Xi_3(x,t) = \frac{k x^{\alpha}}{\Gamma(1+\alpha)} +
\int_0^{\frac{t^{\alpha}}{\Gamma(1+\alpha)}} c(\tau) d\tau$ and
$$
c(\tau)=- k(4k \mu^2(\tau_{\alpha})\diamondsuit S(\tau_{\alpha})-4k
R(\tau_{\alpha})\diamondsuit \mu(\tau_{\alpha}) +6P(\tau_{\alpha})+ k
R^{\diamondsuit 2}(\tau_{\alpha}) /S(\tau_{\alpha}) )/6,
$$
\begin{eqnarray}
V_4(x,t) = -\frac{2k B_2^{\diamondsuit 2}(t)}{S(t)\diamondsuit Q(t)}
+\frac{12k B_{1}(t)\diamondsuit B_{2}^{\diamondsuit 2}(t)}{Q(t)
\diamondsuit S(t) \diamondsuit \left[B_{1}(t)-S(t)\diamondsuit \mu(t)
\diamondsuit \exp^{\diamondsuit}\{
-\frac{B_2(t)}{S(t)}\diamondsuit \Xi_4(x,t) \}\right]}  \mbox{~}\cr
-\frac{12k B_{1}^{\diamondsuit 2}(t) \diamondsuit B_{2}(t)}{Q(t)
\diamondsuit S(t)\diamondsuit \left[B_{1}(t)-S(t)\diamondsuit \mu(t)
\diamondsuit \exp^{\diamondsuit}\{-\frac{B_2(t)}{S(t)}\diamondsuit 
\Xi_4(x,t) \}\right]^{\diamondsuit 2}},
\label{wick-sol-4}
\end{eqnarray}
where $\Xi_4(x,t) = \frac{k x^{\alpha}}{\Gamma(1+\alpha)} +
\int_0^{\frac{t^{\alpha}}{\Gamma(1+\alpha)}} c(\tau) d\tau$ and
$$
c(\tau)=-k(4k \mu^{\diamondsuit 2}(\tau_{\alpha})\diamondsuit
S(\tau_{\alpha})-4k R(\tau_{\alpha})\diamondsuit \mu(\tau_{\alpha})
-P(\tau_{\alpha}) +k R^{\diamondsuit 2}(\tau_{\alpha})/S(\tau_{\alpha})), 
$$
\begin{eqnarray}
V_5(x,t) = \frac{3k B_2^{\diamondsuit 2}(t)}{4S(t)\diamondsuit
Q(t)}-\frac{3k C_1^{\diamondsuit 2}(t) \diamondsuit C_2(t)}{4Q(t)
\diamondsuit S(t) \diamondsuit \left[ C_{1}(t)+S(t) \diamondsuit 
\mu(t)\diamondsuit\exp^{\diamondsuit}\{
-\frac{C_2(t) }{S(t)}\diamondsuit \Xi_5(x,t) \} \right]^{\diamondsuit 2}}, 
\label{wick-sol-5}
\end{eqnarray}
where $\Xi_5(x,t) = \frac{k x^{\alpha}}{\Gamma(1+\alpha)} +
\int_0^{\frac{t^{\alpha}}{\Gamma(1+\alpha)}} c(\tau) d\tau$ and
$$
c(\tau)=-k(12k \mu^{\diamondsuit 2}(\tau_{\alpha})\diamondsuit
S(\tau_{\alpha})-12 k R(\tau_{\alpha})\diamondsuit \mu(\tau_{\alpha}) +8
P(\tau_{\alpha})+3 k R^{\diamondsuit 2}(\tau_{\alpha})/ S(\tau_{\alpha}))/8, 
$$
and
\begin{eqnarray}
V_6(x,t) = -\frac{k B_2^{\diamondsuit 2}(t)}{3S(t)\diamondsuit Q(t)}
+\frac{2k D_1(t) \diamondsuit D_2^{\diamondsuit 2}(t)}{Q(t)\diamondsuit S(t)
\diamondsuit \left[ D_{1}(t)-6S(t) \diamondsuit \mu(t) \diamondsuit 
\exp^{\diamondsuit}
\{-\frac{D_2(t) }{6S(t)} \diamondsuit \Xi_6(x,t) \} \right] }
\mbox{~}\cr -\frac{k D_1^{\diamondsuit 2}(t) \diamondsuit D_2(t)}{3Q(t)
\diamondsuit S(t)
\diamondsuit \left[ D_{1}(t)-6S(t) \diamondsuit \mu(t) \diamondsuit 
\exp^{\diamondsuit}\{
-\frac{D_2(t) }{6S(t)} \diamondsuit \Xi_6(x,t) \} \right]^{\diamondsuit 2}}, 
\label{wick-sol-6}
\end{eqnarray}
where $\Xi_6(x,t) = \frac{k x^{\alpha}}{\Gamma(1+\alpha)} +
\int_0^{\frac{t^{\alpha}}{\Gamma(1+\alpha)}} c(\tau) d\tau$ and
$$
c(\tau)= k(4k \mu^{\diamondsuit 2}(\tau_{\alpha})\diamondsuit
S(\tau_{\alpha})-4k R(\tau_{\alpha})\diamondsuit \mu(\tau_{\alpha}) -6
P(\tau_{\alpha})+ k R^{\diamondsuit 2}(\tau_{\alpha}) /S(\tau_{\alpha}))/6. 
$$
Here, $\tau_{\alpha}=[\tau\Gamma(1+\alpha)]^{1/\alpha}$,
$B_{1}(t)=R(t)-S(t)\diamondsuit \mu(t)$, $B_{2}(t)=R(t)-2S(t)
\diamondsuit \mu(t)$, $C_{1}(t)=R(t)-6S(t) \diamondsuit \mu(t)$, 
$C_{2}(t)=R(t)-5S(t) \diamondsuit \mu(t)$, $D_{1}(t)=R(t)+4S(t)
\diamondsuit \mu(t)$, $D_{2}(t)=R(t)-2S(t)\diamondsuit \mu(t)$, 
and $k$ is a nonzero constant.

Note that the obtained solutions contain arbitrary
functions, which reveal the physical quantities $V_i$, $i=1,2,\ldots,6$. 
These solutions possess rich structures and can be used to
discuss some particular physical situations through the choice 
of the arbitrary functions. Note also that for
different forms of $P(t)$, $Q(t)$, $R(t)$ and $S(t)$, we can get 
white noise functional solutions of exponential type of Eq.
(\ref{rlw-wick}) from the obtained solutions
(\ref{wick-sol-1})--(\ref{wick-sol-6}).

\begin{example} 
Assume $R(t) = 0$ and take $P(t) = f_1 (t)+c_1 W(t)$, 
$Q(t) = f_2 (t)+c_2 W(t)$, $S(t) = f_4 (t)+c_4 W(t)$, 
where $c_i$, $i=1,2,4$, are arbitrary constants 
and $f_i (t)$, $i=1,2,4$, are bounded functions
on $\mathbb{R}^{+}$, where $W(t)$ is the Gaussian white noise, i.e.,
$W(t) = \dot{B}(t) = dB(t)/dt$, $B(t)$ is a Brownian motion. Also we have the
Hermite transform $P(t,z) = f_1 (t)+c_1 \tilde{W}(t,z)$, 
$Q(t,z) = f_2(t)+c_2 \tilde{W}(t,z)$, $S(t,z) = f_4 (t)+c_4 \tilde{W}(t,z)$, 
where $\tilde{W}(t,z) = \sum_{i=1}^{\infty} \int_0^t \zeta_i(s) ds z_i$,
$z = (z_1, z_2, \ldots) \in \mathbb{C}^{n}$ is a parameter
vector and $\zeta_i(s)$ is defined in \cite{holden}. By the
definition of $\tilde{W}(t,z)$, the Wick-type exact traveling wave solutions
(\ref{wick-sol-1})--(\ref{wick-sol-6}) give the white noise functional
solutions of Eq.~(\ref{rlw-wick}) as follows:
\begin{eqnarray}
V_1(x,t) = \frac{48k (f_4(t)+c_4 W(t)) \mu^2(t)}{(f_2(t)+c_2
W(t))\left[ 1+ \exp\{2\mu(t)\Xi_1(x,t) \} \right] } \mbox{~}
-\frac{48k (f_4(t)+c_4 W(t))\mu^2(t)}{(f_2(t)+c_2 W(t)) \left[ 1+
\exp\{2\mu(t)\Xi_1(x,t) \} \right]^2}, \label{ex-1}
\end{eqnarray}
where $\Xi_1(x,t) = \frac{k x^{\alpha}}{\Gamma(1+\alpha)} +
\int_0^{\frac{t^{\alpha}}{\Gamma(1+\alpha)}} c(\tau) d\tau$, 
$c(\tau)=-k(4k \mu^2(\tau_{\alpha}) (f_4(\tau_{\alpha})+c_4 W(t))
+(f_1(\tau_{\alpha})+c_1 W(t)))$ and
$\tau_{\alpha}=[\tau\Gamma(1+\alpha)]^{1/\alpha}$, 
\begin{eqnarray}
V_2(x,t) = -\frac{675k (f_4(t)+c_4 W(t))\mu^2(t)}{(f_2(t)+c_2
W(t))\left[ 6-\exp\{5\mu(t) \Xi_2(x,t)\} \right]^2}, \label{ex-2}
\end{eqnarray}
where $\Xi_2(x,t) = \frac{k x^{\alpha}}{\Gamma(1+\alpha)} +
\int_0^{\frac{t^{\alpha}}{\Gamma(1+\alpha)}} c(\tau) d\tau$,
$c(\tau)=- (3k^2 \mu^2(\tau_{\alpha}) (f_4(\tau_{\alpha})+c_4 W(t))
-2k(f_1(\tau_{\alpha})+c_1 W(t)))/2$ and
$\tau_{\alpha}=[\tau\Gamma(1+\alpha)]^{1/\alpha}$, 
\begin{eqnarray}
V_3(x,t) = \frac{16 k (f_4(t)+c_4 W(t))\mu^2(t)}{(f_2(t)+c_2 W(t))
\left[ 2-3\exp\{\mu(t) \Xi_3(x,t)/3 \} \right] } \mbox{~}
-\frac{16k (f_4(t)+c_4 W(t))\mu^2(t)}{3 (f_2(t)+c_2 W(t))\left[
2-3\exp\{\mu(t) \Xi_3(x,t)/3 \} \right]^2}, \label{ex-3}
\end{eqnarray}
where $\Xi_3(x,t) = \frac{k x^{\alpha}}{\Gamma(1+\alpha)} +
\int_0^{\frac{t^{\alpha}}{\Gamma(1+\alpha)}} c(\tau) d\tau$, 
$c(\tau)=- (4k^2 \mu^2(\tau_{\alpha}) (f_4(\tau_{\alpha})+c_4
W(t))+6k(f_1(\tau_{\alpha})+c_1 W(t)))/6$ and
$\tau_{\alpha}=[\tau\Gamma(1+\alpha)]^{1/\alpha}$, 
\begin{eqnarray}
V_4(x,t) = -\frac{8k \mu^2(t)(f_4(t)+c_4 W(t))}{f_2(t)+c_2 W(t)}
+\frac{48k (f_4(t)+c_4 W(t)) \mu^2(t)}{(f_2(t)+c_2 W(t))\left[ 1+
\exp\{2\mu(t)\Xi_4(x,t) \} \right] } \mbox{~}\cr -\frac{48k
(f_4(t)+c_4 W(t))\mu^2(t)}{(f_2(t)+c_2 W(t)) \left[ 1+
\exp\{2\mu(t)\Xi_4(x,t) \} \right]^2}, \label{ex-4}
\end{eqnarray}
where $\Xi_4(x,t) = \frac{k x^{\alpha}}{\Gamma(1+\alpha)} +
\int_0^{\frac{t^{\alpha}}{\Gamma(1+\alpha)}} c(\tau) d\tau$, 
$c(\tau)=k(4k \mu^2(\tau_{\alpha}) (f_4(\tau_{\alpha})+c_4 W(t))
+(f_1(\tau_{\alpha})+c_1 W(t)))$, and
$\tau_{\alpha}=[\tau\Gamma(1+\alpha)]^{1/\alpha}$, 
\begin{eqnarray}
V_5(x,t) = -\frac{3k \mu^2(t)(f_4(t)+c_4 W(t))}{f_2(t)+c_2 W(t)}
-\frac{675k (f_4(t)+c_4 W(t))\mu^2(t)}{(f_2(t)+c_2 W(t))\left[
6-\exp\{5\mu(t) \Xi_2(x,t)\} \right]^2}, 
\label{ex-5}
\end{eqnarray}
where $\Xi_5(x,t) = \frac{k x^{\alpha}}{\Gamma(1+\alpha)} +
\int_0^{\frac{t^{\alpha}}{\Gamma(1+\alpha)}} c(\tau) d\tau$,
$c(\tau)=- (3k^2 \mu^2(\tau_{\alpha}) (f_4(\tau_{\alpha})+c_4 W(t))
+2k(f_1(\tau_{\alpha})+c_1 W(t)))/2$ and
$\tau_{\alpha}=[\tau\Gamma(1+\alpha)]^{1/\alpha}$, and
\begin{eqnarray}
V_6(x,t) = -\frac{4k \mu^2(t)(f_4(t)+c_4 W(t))}{3(f_2(t)+c_2 W(t))}
+\frac{16 k (f_4(t)+c_4 W(t))\mu^2(t)}{(f_2(t)+c_2 W(t)) \left[
2-3\exp\{\mu(t) \Xi_3(x,t)/3 \} \right] } \mbox{~}\cr -\frac{16k
(f_4(t)+c_4 W(t))\mu^2(t)}{3 (f_2(t)+c_2 W(t))\left[ 2-3\exp\{\mu(t)
\Xi_3(x,t)/3 \} \right]^2}, 
\label{ex-6}
\end{eqnarray}
where $\Xi_3(x,t) = \frac{k x^{\alpha}}{\Gamma(1+\alpha)} +
\int_0^{\frac{t^{\alpha}}{\Gamma(1+\alpha)}} c(\tau) d\tau$,
$c(\tau)=(2k^2 \mu^2(\tau_{\alpha}) (f_4(\tau_{\alpha})+c_4 W(t))
-3k(f_1(\tau_{\alpha})+c_1 W(t)))/2$, and
$\tau_{\alpha}=[\tau\Gamma(1+\alpha)]^{1/\alpha}$. Here $k$ is a
nonzero constant.
\end{example}

\begin{remark}
Let $\alpha=1$, $P(t)=1, Q(t)=1, R(t)=0, S(t)=-1$. Eq.
(\ref{rlw-wick}) becomes the RLW equation \cite{benjamin,99}:
\begin{eqnarray}
\frac{\partial v}{\partial t} + \frac{\partial v}{\partial x} + v
\frac{\partial v}{\partial x} - \frac{\partial^{3} v}{\partial t
\partial x^2} =0. \label{bbm}
\end{eqnarray}
Now, the solutions for the  Eq. (\ref{bbm}) can be obtained 
by the derived solutions of (\ref{rlw-wick}). Here we provide 
some  exact solutions of Eq. (\ref{bbm}) in the following form:
\begin{eqnarray}
v_1(x,t) =  - \frac{48k \mu^2(t)}{\left[ 1+ \exp\{2\mu(t)\xi_1(x,t)\} 
\right]}  +\frac{24k \mu^2(t)}{\left[ 1+ \exp\{2\mu(t)\xi_1(x,t)\} 
\right]^2}, 
\label{rm-1}
\end{eqnarray}
where $\xi_1(x,t) = kx + \int_0^t k(4k \mu^2(\tau)-1) d\tau$, and
\begin{eqnarray}
v_2(x,t) = \frac{135 k \mu^2(t)}{\left[ 6-\exp\{5\mu(t) \xi_2(x,t)\}
\right]^2}, \label{rm-2}
\end{eqnarray}
where $\xi_2(x,t) = k x - \int_0^t \frac{k}{2}(3k \mu^2(\tau)+2)
d\tau$ and
\begin{eqnarray}
v_5(x,t) = -3k\mu^2(t)+\frac{135 k \mu^2(t)}{\left[ 6-\exp\{5\mu(t)
\xi_5(x,t)\} \right]^2}, \label{rm-5}
\end{eqnarray}
where $\xi_5(x,t) = k x + \int_0^t \frac{k}{2}(3k \mu^2(\tau)-2)
d\tau$, here $k$ is a nonzero constant.
\end{remark}

\begin{figure}[ht!]
\begin{center}
\epsfig{file=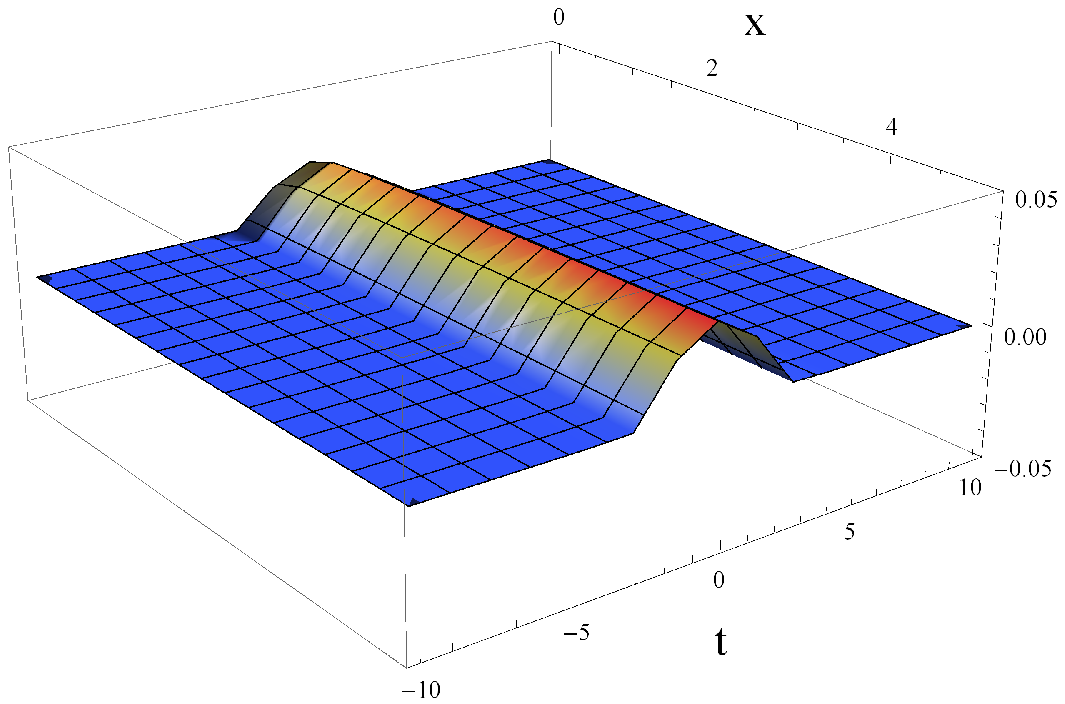, scale=0.44}
\mbox{~}  \epsfig{file=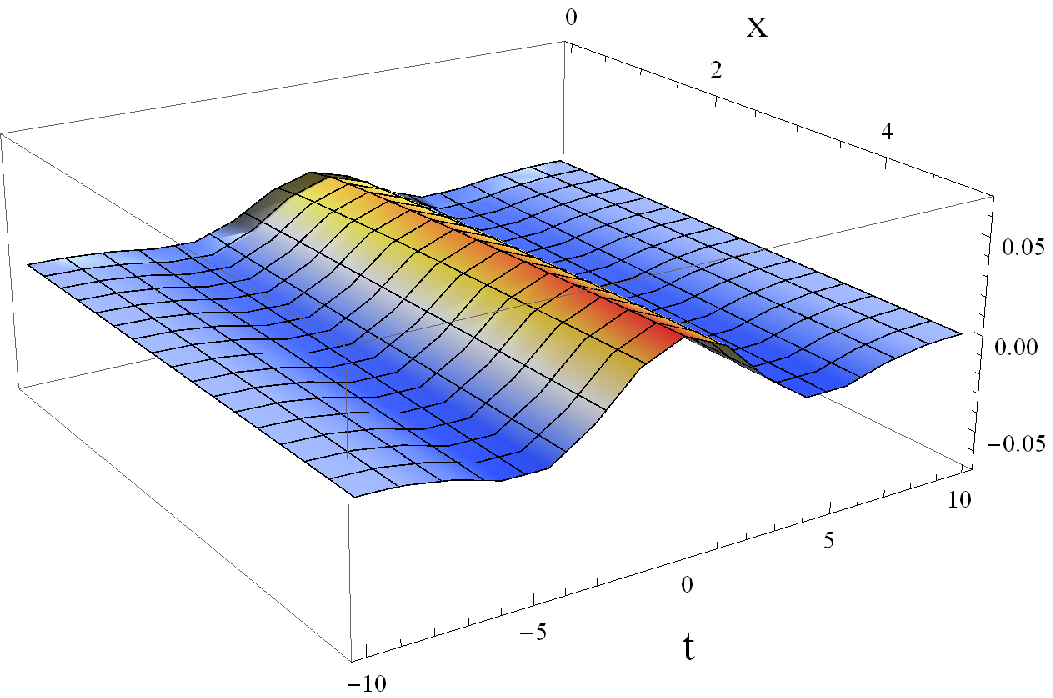, scale=0.44}
\mbox{~}  \epsfig{file=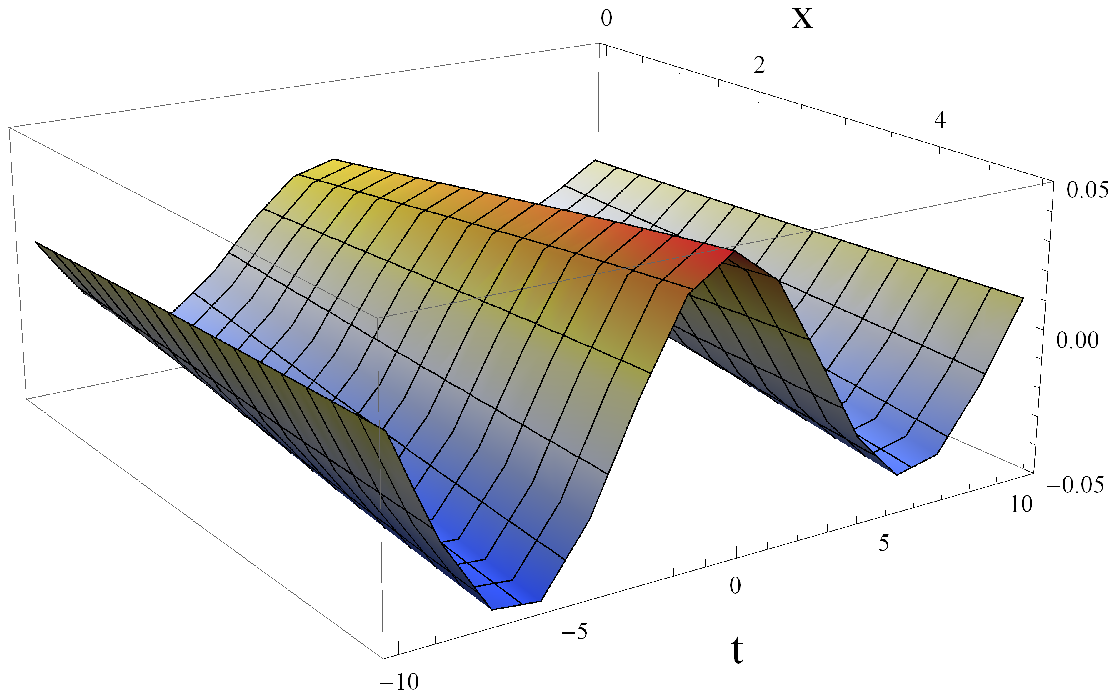, scale=0.44}
\end{center}
\vspace{-0.5cm}
\caption{Wave profiles of exact solution
(\ref{ex-1}) without white noise as $\alpha=0.25, 0.5, 1$.}  
\label{v1-W-0}
\end{figure}

Now, we discuss the solitary waves and the amplitudes of the derived
solutions of Eq. (\ref{rlw-wick}) for fractional-orders with white
noise under some constraints. From the obtained exact solutions 
of Eq. (\ref{rlw-wick}), we perform wave dynamics
related to the damped oscillatory kink wave \cite{300}.
In Fig.~\ref{v1-W-0}, we show the dynamic behavior of exact solution
(\ref{ex-1}) with $W(t)=0$ under $k=0.05$, $\mu(t)=-3$, $f_1(t)=-0.2$,
$f_2(t)=10$, $f_4(t) =0.01 \cos(0.5t)$, $c_1=c_2=c_4=1$ as $\alpha=0.25$,
$0.5$, $1$.
Note that it satisfies  $ V_{1}(x, t) \rightarrow 0$ as $t
\rightarrow \pm \infty$ for all $x$ as $\alpha=0.25, 0.5$ and $
V_{1}(x, t) \rightarrow 0$ as $t \rightarrow \infty$ for all $x$ as
$\alpha=1$.

\begin{figure}[ht!]
\begin{center}
\epsfig{file=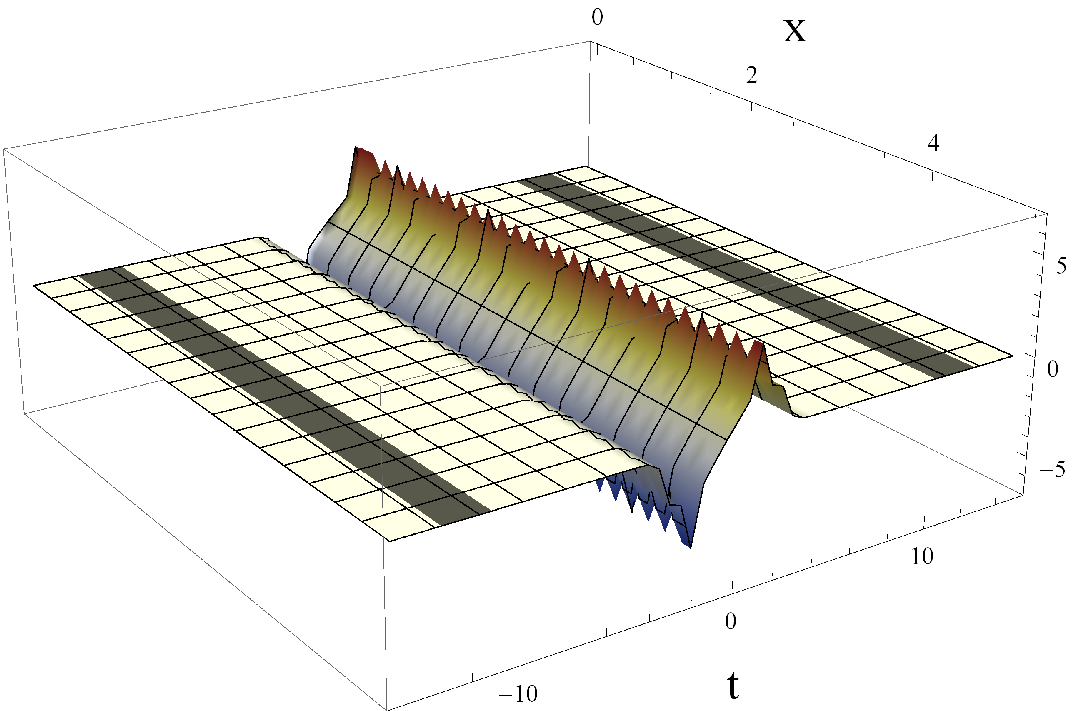, scale=0.45}
\mbox{~}  \epsfig{file=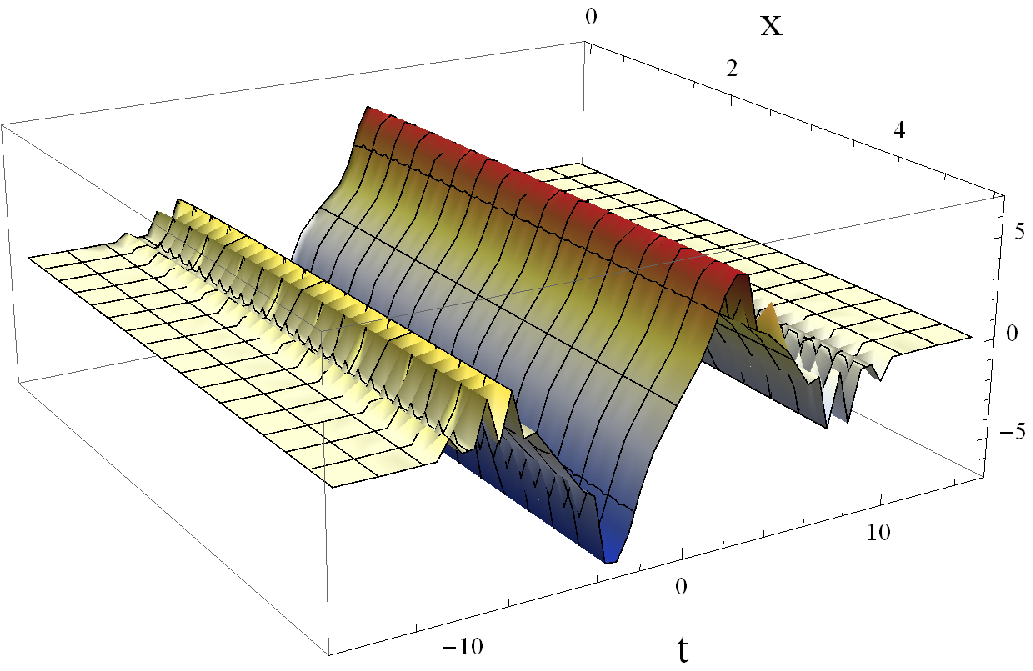, scale=0.45}
\mbox{~}  \epsfig{file=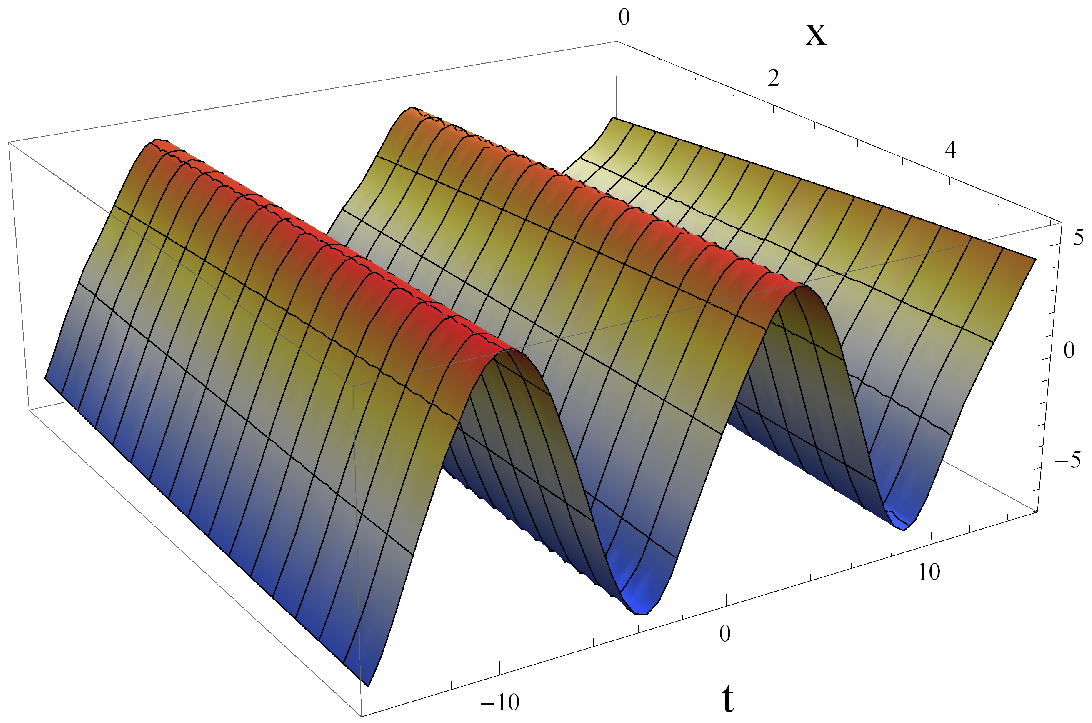, scale=0.45}
\end{center}
\vspace{-0.5cm}
\caption{Wave profiles of exact solution
(\ref{ex-1}) with white noise as $\alpha=0.25, 0.5, 1$.}  
\label{v1-W}
\end{figure}

\begin{figure}[ht!]
\begin{center}
\epsfig{file=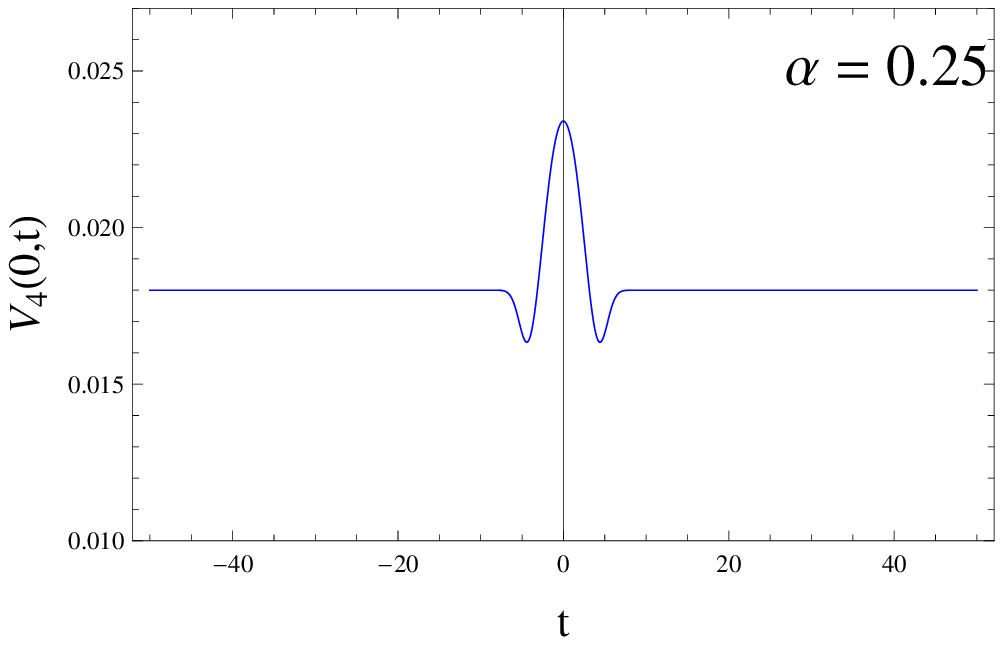, scale=0.48}
\mbox{~}  \epsfig{file=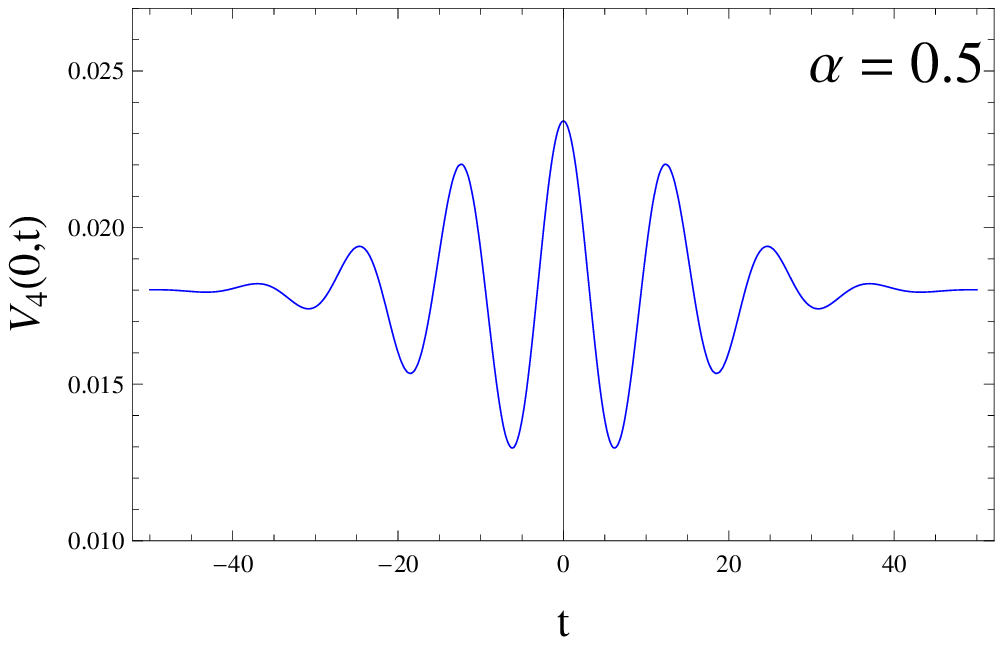, scale=0.48}
\mbox{~} \epsfig{file=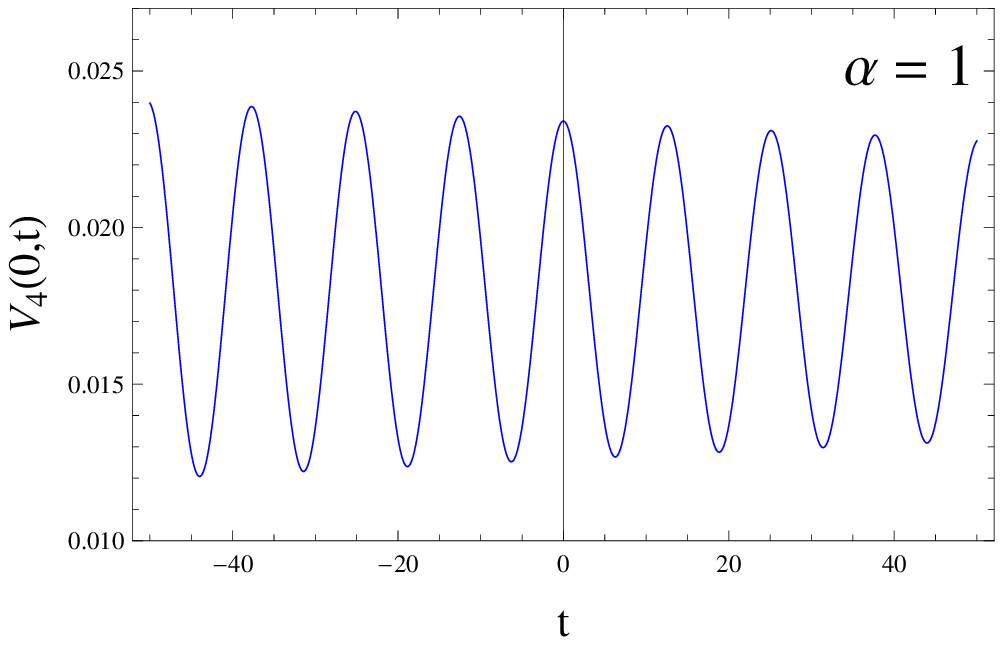, scale=0.48}
\end{center}
\vspace{-0.5cm}
\caption{Wave motions of exact solution (\ref{ex-4}) without 
white noise as $\alpha=0.25, 0.5, 1$.}  
\label{v4-W-0}
\end{figure}

With $k=0.05$, $\mu(t)=-3$, $f_1(t)=-0.2$, $f_2(t)=10$, $f_4(t) =0.01
\cos(0.5t)$, $c_1=c_2=c_4=1$ as $\alpha=0.25$, $0.5$, $1$, $W(t)=
\sin(0.5t)$, we present  the dynamics of the white noise functional
exact solution (\ref{ex-1}) in Fig.~\ref{v1-W}. 
Note that $ V_{1}(x, t) \rightarrow 0$ as $t
\rightarrow \pm \infty$ for all $x$ and the dynamics represent
irregular movements in some interval of $t$ by employing white
noise. For $\alpha=1$, it is like a periodic solution that satisfies
$ V_{1}(x, t) \rightarrow 0$ as $t \rightarrow \infty$ for all $x$. 
In Fig.~\ref{v4-W-0}, we  provided  the dynamics of exact solution
(\ref{ex-4}) with $W(t)=0$, under $k=0.05$, $\mu(t)=-1.5$, $f_1(t)=-0.2$,
$f_2(t)=10$, $f_4(t) =0.01 \cos(0.5t)$, $c_1=c_2=c_4=1$. For
$\alpha=0.25$, $0.5$, $1$, it gives $V_4(x,t) \rightarrow 0.018$ as $t
\rightarrow \infty$ for all $x$.
\begin{figure}[ht!]
\begin{center}
\epsfig{file=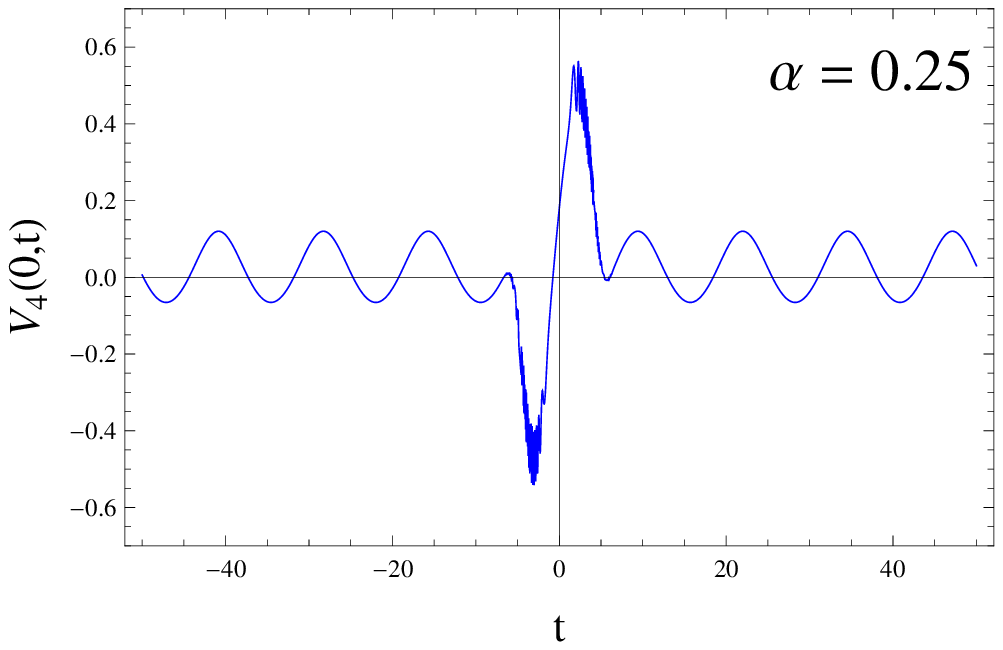, scale=0.48} \mbox{~}
\epsfig{file=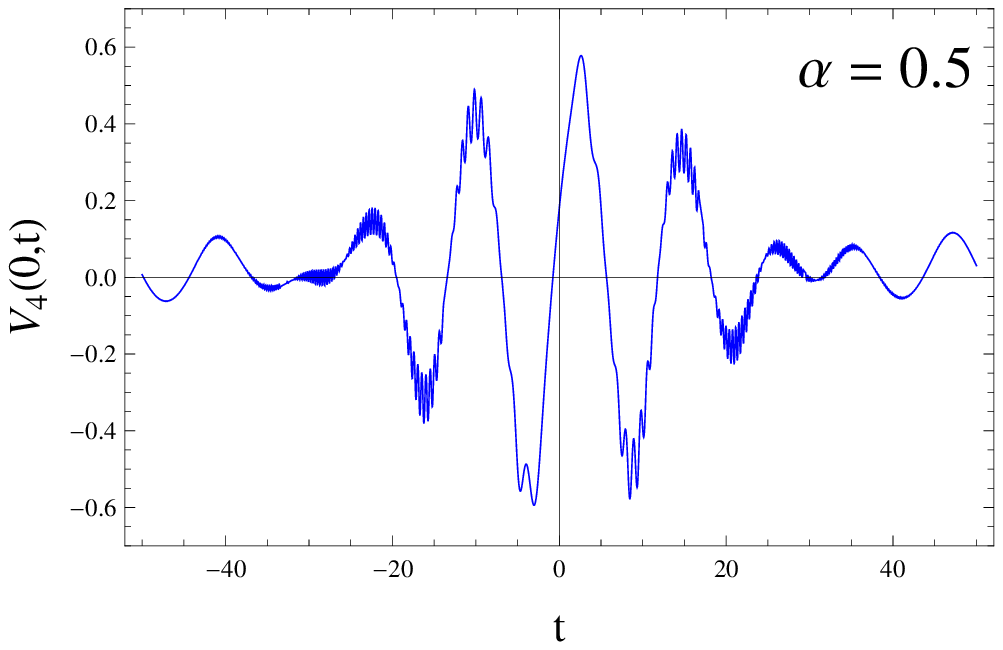, scale=0.48}
\mbox{~} \epsfig{file=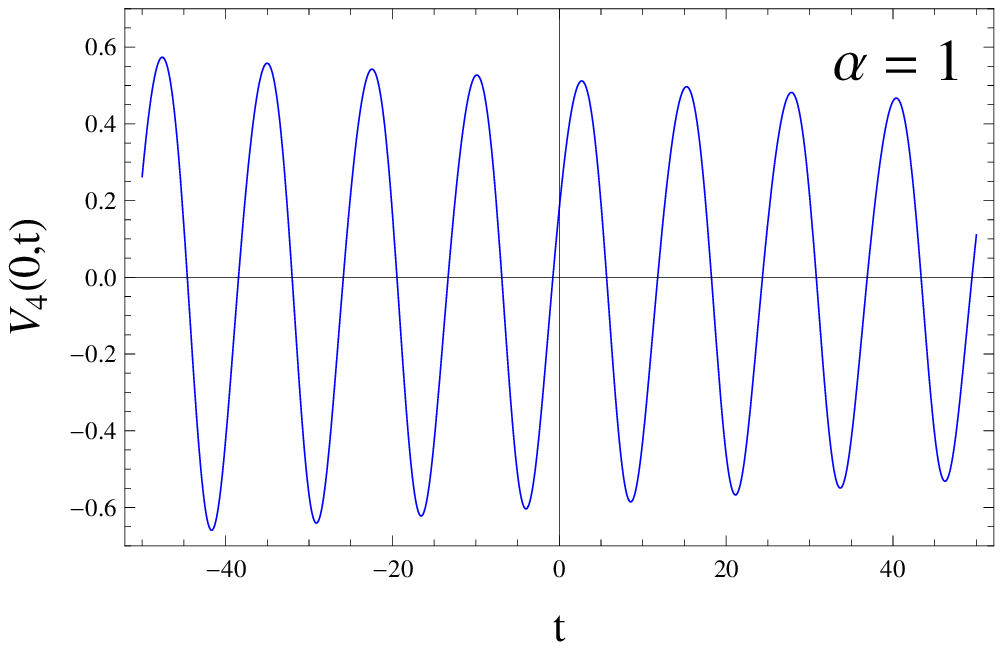, scale=0.48}
\end{center}
\vspace{-0.5cm}
\caption{Wave motions of exact solution (\ref{ex-4}) 
with white noise as $\alpha=0.25, 0.5, 1$.}
\label{v4-W}
\end{figure}
Under a white noise $W(t) = \sin(0.5t)$ and $k=0.05$, $\mu(t)=-1.5$,
$f_1(t)=-0.2$, $f_2(t)=10$, $f_4(t) =0.01 \cos(0.5t)$, $c_1=c_2=c_4=1$, for
$\alpha=0.25, 0.5, 1$, Fig.~\ref{v4-W} gives wave motions of the
white noise functional solution of (\ref{ex-4}). It can be seen that the
behaviors of exact solution $V_4(x,t)$ are roughly dynamic in
$(-t_0, t_0)$ and it gives periodic wave motions in $(\infty, -t_0)$
and $(t_0, \infty)$ for $\alpha=0.25$ and $\alpha=0.5$. Furthermore,
for $\alpha=1$, the values of $V_4(0,t)$ lies between $-0.08$ and
$0.12$ as $t$ approaches $\infty$. From the simulations, 
it is concluded that we obtain different dynamics of the white noise 
functional solutions of Eq. (\ref{rlw-wick}) by the role of
white noise as fractional orders.


\section{Conclusion}
\label{sec:4}

In this paper, we obtained Wick-type exact solutions 
of the stochastic nonlinear Schr\"{o}dinger 
and the Wick-type stochastic fractional RLW-Burgers equations 
by employing an improved computational method with
integral and fractional order. Specifically, the Hermite
transform is implemented, which changes Wick products into ordinary
products and consequently converts Wick-type stochastic
equations to deterministic nonlinear ordinary differential ones. 
Then we have applied the improved computational technique 
for determining exact traveling wave solutions for associated 
deterministic equations. Finally, the inverse Hermite transform 
is applied to obtained solutions for obtaining 
Wick-type exact solutions in the Wick type stochastic case.


\section*{Acknowledgements} 

This research was supported by the Basic Science Research Program 
of the National Research Foundation of Korea (NRF), funded by the 
Ministry of Education, no. NRF-2019R1A6A1A10073079. 
Debbouche and Torres were supported by 
FCT within project UID/MAT/04106/2019 (CIDMA). 
 


\end{document}